\documentclass[11pt]{article}

\usepackage{amsmath}
\usepackage{amssymb}
\usepackage{amscd}
\usepackage{palatino}
 \usepackage{euler}
\usepackage{latexsym}
\usepackage{theorem}
\usepackage{setspace}
\usepackage{epsfig}
\usepackage{graphics}
\usepackage{psfrag}
\usepackage{amsfonts}
\usepackage{nomencl}

\newtheorem{theorem}{Theorem}[section]
\newtheorem{lemma}[theorem]{Lemma}
\newtheorem{proposition}[theorem]{Proposition}
\newtheorem{corollary}[theorem]{Corollary}

{\theorembodyfont{\rmfamily}
\theoremstyle{plain}
\newtheorem{definition}[theorem]{Definition}
\newtheorem{example}[theorem]{Example}

\newtheorem{remark}[theorem]{Remark}
}

\numberwithin{equation}{section}
\numberwithin{figure}{section}

\input xy
 \xyoption{all}

\usepackage{multicol}

\makeatletter

\def\theglossary{\@restonecoltrue\if@twocolumn\@restonecolfalse\fi
\columnseprule\z@ \columnsep 35\p@
\let\@makessectionhead\indexsec
\@xp\section\@xp*\@xp{\glossaryname}%
\let\item\@idxitem
\parindent\z@  \parskip\z@\@plus.3\p@\relax
\footnotesize}

\def\glossaryname{Notation Index}

\makeatother

\makeglossary

\newcommand{\liet}{{\mathfrak t}}

\newcommand{\e}{\varepsilon}

\newcommand{\codim}{\operatorname{codim}}

\renewcommand{\H}{{\mathbb{H}}}
\newcommand{\C}{{\mathbb{C}}}
\newcommand{\Zt}{{\mathbb{Z}_2}}
\newcommand{\Z}{{\mathbb{Z}}}
\newcommand{\Q}{{\mathbb{Q}}}
\newcommand{\R}{{\mathbb{R}}}

\newcommand{\Hn}{{\H^n}}
\newcommand{\Cn}{{\C^n}}

\renewcommand{\cot}{T^*\Cn}
\newcommand{\Tk}{T^k}
\newcommand{\Tn}{T^n}

\newcommand{\Td}{T^d}

\renewcommand{\t}{\mathfrak{t}}       
\newcommand{\algt}{\mathfrak{t}}       
\newcommand{\tk}{\mathfrak{t}^k}
\newcommand{\tn}{\mathfrak{t}^n}
\newcommand{\td}{\mathfrak{t}^d}

\newcommand{\tnd}{(\tn)^*}
\newcommand{\tdd}{(\td)^*}

\newcommand{\subs}{\subseteq}

\newcommand{\into}{{\hookrightarrow}}

\newcommand{\htmzt}{H^{2*}_T(M;\Zt)}

\newcommand{\htrx}{H^*_{T_\R}(Q;\Zt)}

\newcommand{\hs}{\hspace{3pt}}

\renewcommand{\mod}{{/\!\!/}}

\newcommand{\mc}{{\mathcal{C}}}
\newcommand{\arr}{{\mathcal{H}}}

\newcommand{\iso}{\cong}

\newcommand{\<}{\left<}
\renewcommand{\>}{\right>}

\newcommand{\mr}{\mu_{\R}}
\renewcommand{\mc}{\mu_{\C}}

\newcommand{\qed}{\hfill \mbox{$\Box$}\medskip\newline}

\newenvironment{proof}{\noindent {\bf Proof:}}{\qed \par}

\begin{document}
\begin{spacing}{1.1}

\noindent
{\LARGE \bf The equivariant cohomology of hypertoric varieties and
their real loci}
\bigskip\\
{\bf Megumi Harada}\footnote{{\tt megumi@math.toronto.edu}}\\
Department of Mathematics, University of Toronto,
Toronto, Ontario, M5S 3G3 Canada\smallskip \\
{\bf Tara S. Holm }\footnote{{\tt tsh@math.berkeley.edu}.
\newline \mbox{~~~~} The second author was supported by a National
Science Foundation Postdoctoral Research Fellowship.
\newline \mbox{~~~~} {\it MSC 2000 Subject Classification}:
Primary: 55N91 \hspace{0.1in} Secondary: 53C26, 05C90
\newline \mbox{~~~~} {\it Keywords:} equivariant cohomology,
symplectic, real loci, toric varieties, hyperk\"ahler} \\
Department of Mathematics, University of California,
Berkeley, CA 94720, USA
\bigskip

{\small
\begin{quote}
\noindent {\em Abstract.} Let $M$ be a Hamiltonian $T$ space
with a proper moment map, bounded below in some component.
In this setting, we give a combinatorial description of the
$T$-equivariant cohomology of $M$, extending results of Goresky,
Kottwitz and MacPherson and techniques of Tolman and
Weitsman. Moreover, when $M$ is equipped with an antisymplectic
involution $\sigma$ anticommuting with the action of $T$, we also
extend to this noncompact setting the
``mod 2'' versions of these results to the real locus $Q:= M^\sigma$
of $M$. We give applications
of these results to the theory of hypertoric varieties.
\end{quote}
}
\bigskip

\tableofcontents

\section{Introduction}\label{sec:intro}

In this article, we present two main results and demonstrate their use
through several explicit computations. In the first result, we generalize
to the noncompact setting a theorem of Goresky, Kottwitz, and
MacPherson that computes $T = T^n$-equivariant cohomology rings of
compact Hamiltonian $T$ spaces satisfying some technical
conditions \cite{GKM}.  Suppose in addition that $M$ is equipped with
an antisymplectic involution $\sigma$ that anticommutes with the
$T$ action.  In the second result, we generalize to the noncompact
setting theorems \cite{BGH, Dui83, GolHol04, Sch} that
compute the equivariant cohomology of the real locus $Q:=M^\sigma$ of a
Hamiltonian $T$ space $M$ satisfying similar technical conditions. The
motivating examples for this paper are the hypertoric
varieties studied in \cite{BieDan, Kon99, Kon00,
HauStu, HP02} and their real loci.  We present these in detail.

We first recall the basic idea of the theorem of Goresky,
Kottwitz, and MacPherson, which we call the GKM theorem. For a
compact Hamiltonian $T$ space $M$,
Kirwan showed \cite{Kir} that the inclusion $M^T\hookrightarrow M$
induces an injection $H_T^*(M;\Q)\hookrightarrow H_T^*(M^T;\Q)$ in
equivariant cohomology.
Since $T$ acts trivially on $M^T$, when $M^T$ consists of isolated
points,  the ring
\(H^*_T(M^T;\Q)\) is a direct product of polynomial rings
\[H^*_T(M^T;\Q) \iso \prod_{p \in M^T} H^*_T(p;\Q) \iso \prod_{p
\in M^T} Sym(\algt^*).\] Hence, in order to compute the
equivariant cohomology ring $H^*_T(M;\Q)$, it suffices to identify
its image in $H^*_T(M^T;\Q)$.  Suppose in addition that the $T$-isotropy weights
$\{\alpha_{p,i}\}$ are pairwise linearly independent at each fixed
point $p$. The GKM theorem then asserts that the image of $H^*_T(M;\Q)$
in $H^*_T(M^T;\Q)$ is the same as that of the one-skeleton of $M$, which in
turn can be described combinatorially in terms of a graph $\Gamma$
and the $T$-isotropy data. Thus, the
computation of $H^*_T(M;\Q)$ is translated into a problem of
combinatorics.

In the setting of compact Hamiltonian $T$ manifolds equipped with
an additional antisymplectic involution $\sigma$, we define the
{\em real locus} $Q := M^{\sigma}$ of $M$, which is equipped with an
action of the discrete subgroup $T_{\R} := (\Zt)^n$ of $T$.
The mod $2$ GKM
theorem \cite{BGH, Sch} gives a similar  combinatorial description of the
image of the equivariant cohomology of the real locus $Q :=
M^{\sigma}$ as a subring of the equivariant cohomology (with $\Zt$
coefficients) of $Q^{T_\R}$.

Examples of Hamiltonian $T$ spaces satisfying the GKM hypotheses
include coadjoint orbits of compact Lie groups and toric varieties.
In the case of coadjoint orbits, the combinatorial description given
by the GKM theorem has proved useful in the theory of equivariant
Schubert calculus (for example, see \cite{KnuTao03}).
These examples also have natural antisymplectic involutions, and the
mod $2$ results apply to these examples.   The combinatorial
description of the equivariant cohomology of real loci of certain
toric varieties has applications to string theory (see
\cite{BGH}).

Thus far, we have required $M$ to be compact.  However, there are many
noncompact examples that nonetheless fit into this framework.  For
example, hypertoric varieties \cite{BieDan, Kon99, Kon00, HauStu,
HP02} equipped with a $T\times S^1$ action exhibit many of the
properties of 
compact K\"ahler toric varieties. For instance, like their K\"ahler
counterparts, the $T\times S^1$-isotropy weights at each fixed point
are pairwise linearly independent. The hypertoric varieties are also
equipped with a natural antisymplectic involution, and the computation
of the equivariant cohomology of the real loci has applications to the
theory of hyperplane arrangements \cite[Section 5]{HP02}. Moreover,
coadjoint orbits of affine Kac-Moody algebras equipped with an
appropriate $T\times S^1$ action satisfy pairwise linear independence
properties.  Thus, the examples strongly motivate us to demonstrate
GKM and mod 2 GKM theorems in the noncompact setting.

The essential observation in this paper is that the Morse theoretic
arguments given in \cite{TW99} work with only slight modifications in
the setting when there is a direction of the moment map which is
proper and bounded below. These two hypotheses allow us to apply the
same local Morse theoretic arguments: the properness ensures the
compactness of critical sets, and the boundedness allows us to make an
inductive argument by providing a base case.  This is the case for
some of the examples mentioned above; for smooth
hypertoric varieties, it is indeed true that there exists a component
of the $T\times S^1$ moment map which is proper and bounded below
\cite{HP02}.  Tor the coadjoint orbits of affine Kac-Moody
algebras, however, the methods presented in this paper do not suffice. Different
arguments must be used to give a GKM description for these examples
\cite{HHH04}.  The results in \cite{HHH04} are phrased in the language
of cell complexes, but they also achieve a combinatorial description
of equivariant cohomology. We also note that the results in this paper
are stated over $\Z$ instead of $\Q$. This changes the statements of
some of the technical hypotheses on the $T$-isotropy weights.

We now give a brief outline of the contents of this paper. In
Section~\ref{sec:gkm}, we state and prove a GKM theorem in the setting
of noncompact spaces in Theorem~\ref{thm:GKMnoncompact}. We use this
theory in Section~\ref{sec:hypertoric} to analyze in detail the
example of smooth hypertoric varieties equipped with a Hamiltonian
$T^d \times S^1$ action. In particular, we give an isomorphism between
the quotient description of the $T^d \times S^1$-equivariant
cohomology of a hypertoric variety given in \cite{HP02} with the GKM
description in Theorem~\ref{thm:quotient_GKM_isom}, and compute
several examples. Further, although the $T^d$ action on $M$ does not
satisfy the GKM hypotheses, we use a ``GKM in stages'' argument to give a
computation of the $T^d$-equivariant cohomology of $M$ by using our
GKM description of its $T^d \times S^1$-equivariant cohomology. 
In Section~\ref{sec:mod2gkm}, we state and prove a
mod $2$ GKM theorem in the noncompact setting, stated in
Theorem~\ref{thm:mod2gkm}. We use this to analyze the real locus of
hypertoric varieties in Section~\ref{sec:realex}. In particular, we
explicitly identify the isomorphism between the $T^d \times
S^1$-equivariant cohomology of a hypertoric variety and the $T_\R^d
\times \Zt$-equivariant cohomology of its real locus in
Proposition~\ref{prop:htRealLocusIsom}. We also mention an application
of these results that is used in \cite{HP02}.

\section{GKM theory for noncompact spaces}\label{sec:gkm}

The goal of this section is to extend results about the
equivariant topology of compact symplectic Hamiltonian manifolds
to situations where the manifold is not necessarily compact.  We replace the compactness hypothesis by
a hypothesis on the moment map: we require it to be proper and
bounded in some direction.  This hypothesis ensures that we can
still use components of the moment map to study the Hamiltonian
manifold Morse theoretically.

Our proofs of Theorems~\ref{thm:injectivity} and~\ref{th:oneskeleton}
follow the outline of the arguments given in \cite{TW99}. The
technical heart of the argument is a lemma due to Atiyah and Bott.
The hypothesis on the moment map ensures that this lemma still
applies to our noncompact setting.  We use this lemma, along
with the Morse theory of the moment map, to show that the
equivariant cohomology of $M$ injects into the equivariant
cohomology of the fixed point set.  We then show that the image is
the same as the image of the equivariant cohomology of the
one-skeleton, $\overline{N}$.  The main result then follows as a
corollary to this: we give the combinatorial description of
$H_T^*(M)$, given additional hypotheses on $M^T$ and on
$\overline{N}$.

We first present the key lemma of Atiyah and Bott.  It is stated
in \cite[Proposition~13.4]{AB82},
\cite[Proposition~5.3.7]{AllPup}.

\begin{lemma}[Atiyah-Bott]\label{lem:AB}
Let ${\cal E}\to B$%
\glossary{E@${\cal E}$, total space of a vector
  bundle}%
\glossary{B@$B$, base space of a vector bundle}
be a complex rank $\ell$ vector bundle over a compact oriented manifold
$B$. Let $T$\glossary{T@$T$, compact $n$-dimensional torus}
be the compact torus $T = (S^1)^d$.
Suppose that $T$ acts on ${\cal E}$ with fixed point set
precisely $B$.  Suppose further that the cohomology of $B$ has no
torsion over $\Z$.  Choose a $T$-invariant Riemannian metric on
${\cal E}$, and let $D$ and $S$ be the corresponding disk and sphere bundles,
respectively, of ${\cal E}$.  Then the long exact sequence of the pair
$(D,S)$ splits into short exact sequences
$$
\xymatrix{ 0 \ar[r] &  H_{T}^*(D,S; \Z) \ar[r] &  H_{T}^*(D;\Z) \ar[r] &
H_{T}^*(S;\Z) \ar[r] & 0 }.
$$
\end{lemma}

\begin{remark}
An alternative statement of this lemma is that
the $T$-equivariant Euler class of the bundle ${\cal E}$ is not a
zero divisor.
\end{remark}

\begin{remark}
The hypothesis that the cohomology of $B$ has no torsion over $\Z$ can
be relaxed to the hypothesis that it has
no $2$-torsion when we take the coefficient ring to be $\Zt$, and
can be removed entirely if we take coefficient ring $\Q$.
\end{remark}

We now turn our attention to finite-dimensional Hamiltonian
$T$ spaces. Suppose that a torus $T$ acts on a symplectic manifold
$M$ in a Hamiltonian fashion.  Then components of the moment map
$\mu : M\to \algt^*$ are Morse-Bott functions on $M$, and if the
component is generic, the critical set is precisely the fixed
point set.  When we assume that a generic component is proper, 
then the connected components of the fixed point set are
compact.  Thus, we may use Lemma~\ref{lem:AB} to study the normal
bundles to these fixed point sets to prove the following
proposition.

\begin{proposition}\label{prop:key}
Let a torus $T$ act on a symplectic manifold $M$ with moment map $\mu
: M \to \algt^*$\glossary{mu@$\mu$, moment map} that is proper in some generic direction \(f := \mu^{\xi}.\)%
\glossary{xi@$\xi$, generic direction in $\t$} \glossary{f@$f$, a
generic component of the moment map $\mu$} Let $c$\glossary{c@$c$,
a critical value of moment map \(f=\mu^{\xi}\)} be a critical
value of $f$. Let $\Sigma_c$\glossary{Sigmac@$\Sigma_c$, component
of the fixed point set $\Sigma=M^T$} be the component of
$\Sigma:=M^T$\glossary{Sigma@$\Sigma$, fixed point set of $M$
under action of $T$} with $\mu^\xi (\Sigma_c ) = c$, and assume
that the cohomology of $\Sigma_c$ has no torsion over $\Z$. For
small $\varepsilon > 0$, let $M_c^\pm := f^{-1}(-\infty , c\pm
\varepsilon)$\glossary{Mcpm@$M_c^\pm$, the inverse image
$f^{-1}(-\infty, c\pm \varepsilon)$}. Then the long exact sequence
of the pair $(M_c^+,M_c^-)$ splits into short exact sequences
$$
\xymatrix{ 0\ar[r] & H_T^*(M_c^+,M_c^-;\Z) \ar[r] & H_T^*(M_c^+;
\Z) \ar[r]^{k^*} & H^*_T(M_c^-;\Z) \ar[r] & 0 }.
$$
Moreover, the restriction from $H_T^*(M_c^+; \Z)$ to $H_T^*(\Sigma_c;
\Z)$ induces an isomorphism from the kernel of $k^*$ to the
classes of $H_T^*(\Sigma_c; \Z)$ that are multiples of
$\tau_c$\glossary{tauc@$\tau_c$, equivariant Euler class of negative
  normal bundle at the critical value $c$} the
equivariant Euler class of the negative normal bundle to $\Sigma_c$.
\end{proposition}

\begin{proof}
This argument appears in \cite{TW99}.  Let $D_c$ and $S_c$ denote the
disc and sphere bundles of the negative normal bundle to the fixed set
$\Sigma_c$.  Using the retraction of the pair $(M_c^+,M_c^-)$ to the
pair $(D_c,S_c)$ and the Thom isomorphism we get the
commutative diagram
$$
\begin{array}{c}
\xymatrix{ {}\ar[r] &
H^{*}_{T}(M_c^+,M_c^-;\Z)\ar[r]\ar[d]^{\iso} &
H^*_{T}(M_c^+;\Z) \ar[r]^{k^*}\ar[d] & H^*_{T}(M_c^-;\Z)\ar[r] &  \\
&  H_{T}^*(D_c,S_c;\Z)\ar[r]\ar[d]_{\mathrm{Thom\ Iso}}^{\iso} & H_{T}^*(D_c;\Z) & & \\
&  H^{*-\lambda}_{T}(D_c;\Z) \ar[ur]_{\cup \tau_c} & & & }
\end{array}.
$$
By the Atiyah-Bott Lemma, the cup product with $\tau_c$ is injective;
therefore the top long exact sequence splits into short exact
sequences.  The proposition follows by a diagram chase.
\end{proof}

Using this proposition, we prove by an inductive argument that the
equivariant cohomology of $M$ injects into the equivariant
cohomology of the fixed point set $\Sigma$.  In order to start the
induction, we now add the hypothesis that a generic component of
the moment map is not only proper, but also bounded below.

\begin{theorem}\label{thm:injectivity}
Let a torus $T$ act on a symplectic manifold $M$ with moment map
$\mu : M\to \algt^*$ that is proper and bounded below in some
generic direction. Suppose that $\Sigma = M^T$ has only finitely
many connected components.  Let $\imath : \Sigma\to M$
\glossary{incl@$\imath$, the inclusion of the fixed point set} be
the inclusion of the fixed point set into $M$. Then the pullback
map
$$
\imath^* : H_T^*(M;\Z) \to H_T^*(\Sigma;\Z)
$$
is injective.
\end{theorem}

\begin{proof}
Choose $\xi$ with $f=\mu^\xi$ generic, proper, and bounded below.
The critical sets of $f$ are precisely the connected components of
$\Sigma$. Thus, by assumption on $\Sigma$, there are only finitely
many critical values of $f$.  Order these critical values as $c_1
< c_2 < \cdots < c_m$, and let $\Sigma_{c_1}, \dots ,
\Sigma_{c_m}$ denote the corresponding critical submanifolds.
These critical submanifolds are compact, since $f$ is proper.  Let
$\Sigma_{c_i}^\pm:=M_{c_i}^\pm\cap \Sigma$. We now proceed by
induction on the critical values.

Let $\varepsilon>0$ be smaller than any of the values
$c_i-c_{i-1}$. The injectivity result is true for $M_{c_1}^+$, as
it is equivariantly homotopic to $\Sigma_{c_1}$. Now assume by
induction that it is true for $M_{c_i}^-$. We will prove that it
is true for $M_{c_i}^+$. Note that $M_{c_i}^-$ is homotopy
equivalent to $M_{c_{i-1}}^+$. We have the long exact sequence of
the pair $(M_{c_i}^+,M_{c_i}^-)$, but by
Proposition~\ref{prop:key}, this splits into short exact
sequences. Thus, we have a commutative diagram
\begin{equation}\label{eq:commdiagr}
\begin{array}{c}
\xymatrix{
0  \ar[r] & H_T^*(M_{c_i}^+,M_{c_i}^-;\Z) \ar[r]\ar[d]
& H_T^*(M_{c_i}^+;\Z) \ar[r]\ar[d]_{\imath^*_+} & H_T^*(M_{c_i}^-;\Z)
\ar[r]\ar[d]_{\imath^*_-}
& 0 \\
0  \ar[r] & H_T^*(\Sigma_{c_i};\Z)\ar[r] &
H_T^*(\Sigma_{c_i}^+;\Z) \ar[r] & H_T^*(\Sigma_{c_i}^-;\Z) \ar[r]
& 0 }
\end{array},
\end{equation}
where we identify $H_T^*(\Sigma_{c_i}^+,\Sigma_{c_i}^-;\Z)\iso
H_T^*(\Sigma_{c_i};\Z)$. The left vertical arrow is an injection,
induced by the Thom isomorphism, and the right vertical arrow is
an injection by the induction hypothesis. A simple diagram chase
shows that the middle vertical arrow must also be an injection.
Since there are only finitely many critical values, the result now
follows by induction.
\end{proof}

Since $T$ is acting trivially on $\Sigma$, the ring
$H_T^*(\Sigma;\Z)$ is isomorphic to the ring
$H^*(\Sigma;\Z)\otimes H_T^*(pt;\Z)$. In general,
$H_T^*(\Sigma;\Z)$ is easier to compute than $H_T^*(M;\Z)$. Thus,
in order to compute $H_T^*(M;\Z)$ as a ring, it now suffices to
describe the image in $H_T^*(\Sigma;\Z)$. We will now show that in
fact the image of $H_T^*(M;\Z)$ is the same as the image of the
equivariant cohomology of a certain subset of $M$.

Let $N$ denote the subset of $M$ given by
$$
N := \{x \in M\ |\ \codim(Stab(x)) = 1\}.
$$ Thus $N$\glossary{N@$N$, the interior of the one-skeleton of $M$}
consists of the points in $M$ whose $T$ orbit is exactly
one-dimensional. We now define the {\em one-skeleton} of $M$ to be
the closure of $N$.  That is, it is the set
$$
\overline{N} := \{ x\in M\ |\ \codim(Stab(x))\leq 1\}.
$$
\glossary{Ncl@$\overline{N}$, the one-skeleton of $M$}
We have the diagram of inclusions
\begin{equation}\label{eq:oneskel}
\begin{array}{c}
\xymatrix{ \Sigma \ar@{^{(}->}[rr]^{\imath}
\ar@{^{(}->}[dr]_{\jmath}\glossary{jmath@$\jmath$, inclusion of
  one-skeleton into $M$} & & M \\
 & \overline{N}\ar@{^{(}->}[ur] & \\
 }
 \end{array}.
\end{equation}
The next theorem states that the image of $H_T^*(M;\Z)$ is the
same as the image of $H_T^*(\overline{N};\Z)$ in
$H_T^*(\Sigma;\Z)$. It is a noncompact version of a theorem of
Tolman and Weitsman \cite[Theorem 1]{TW99}.  As above, the
compactness hypothesis is replaced by the hypothesis that some
generic component of the moment map be proper and bounded. Note
that our theorem holds with $\Z$ coefficients in contrast to
\cite[Theorem 1]{TW99}, which is stated for $\Q$ coefficients.  To
achieve this, we have added an assumption on the $T$ weights on
the negative normal bundle.

\begin{theorem}\label{th:oneskeleton}
Let a torus $T$ act on a symplectic manifold $M$ with moment map $\mu
: M \to \algt^*$ that is proper and bounded below in some generic
direction \(f := \mu^{\xi}.\) Suppose that $\Sigma := M^T$ has only
finitely many connected components.  Suppose further that the distinct
weights of the $T$ action on the negative normal bundle, with respect
to $f$, to any fixed point component are pairwise relatively prime in
$H_T^*(pt;\Z)$.  Then in the diagram in equivariant cohomology,
induced by the inclusions \eqref{eq:oneskel},
$$
\begin{array}{c}
\xymatrix{
H_T^*(M;\Z) \ar[rr]^{\imath^*} \ar[dr] & & H_T^*(\Sigma;\Z) \\
 & H_T^*(\overline{N};\Z) \ar[ur]_{\jmath^*} &
 }
 \end{array},
$$
the image of $\imath^*$ is equal to the image of $\jmath^*$ in
$H_T^*(\Sigma;\Z)$.
\end{theorem}

\begin{remark}
Note that if $T=S^1$, we have $\overline{N} = M$, and the theorem
automatically holds.
\end{remark}

\begin{proof}
We proceed by
induction on the critical values $c_1 < c_2 < \cdots < c_m$ of
$f:=\mu^\xi$.  We first set up our notation. Let $c$ be one of the
critical values of $f$. Define the sets $\overline{N}_c^\pm :=
\overline{N}\cap M_c^\pm$.  Then we have inclusions
$$
\begin{array}{c}
 \xymatrix{
\Sigma_c^\pm \ar@{^{(}->}[rr]^{\imath_\pm} \ar@{^{(}->}[dr]_{\jmath_\pm} & & M_c^\pm \\
 & \overline{N}_c^\pm \ar@{^{(}->}[ur] &
 }\end{array},
$$
which induce, in equivariant cohomology,
$$
\begin{array}{c}
\xymatrix{
H_T^*(M_c^\pm;\Z) \ar[rr]^{\imath_\pm^*} \ar[dr] & & H_T^*(\Sigma_c^\pm;\Z) \\
 & H_T^*(\overline{N}_c^\pm;\Z) \ar[ur]_{\jmath_\pm^*} &
 }\end{array}.
$$
The base case consists of analyzing these diagrams for the minimum
critical value $c_1$.  In this case, $M_{c_1}^-$ and $N_{c_1}^-$
are empty, and $M_{c_1}^+$ and $N_{c_1}^+$ both equivariantly
retract onto $\Sigma_{c_1}$. Thus, both $\imath_+^*$ and
$\jmath_+^*$ are isomorphisms, and therefore have the same image.

Assume now by induction that the statement holds for
$M_{c_{i-1}}^+$. Let $r$ denote the natural restriction from
$im(\jmath_+^*)\subseteq H_T^*(\Sigma_{c_i}^+;\Z)$ to
$H_T^*(\Sigma_{c_i}^-;\Z)$.  Note that the image of $r$ is
contained in $im(\jmath_-^*)$.  By abuse of notation, we will let
$\ker(r)$ denote the inverse image inside
$H_T^*(\Sigma_{c_i};\Z)\iso
H_T^*(\Sigma_{c_i}^+,\Sigma_{c_i}^-;\Z)$ of the kernel of $r$,
using the short exact sequence of the pair
$(\Sigma_{c_i}^+,\Sigma_{c_i}^-)$. Thus, we have a commutative
diagram
\begin{equation}\label{eq:inductiveStep}
 \begin{array}{c}
 \xymatrix{
0  \ar[r] & H_T^*(M_{c_i}^+,M_{c_i}^-;\Z) \ar[r]\ar@{..>}[d]  &
H_T^*(M_{c_i}^+;\Z) \ar[r]\ar[d]_{\imath^*_+} &
H_T^*(M_{c_i}^-;\Z) \ar[r]\ar[d]_{\imath^*_-}
& 0 \\
0  \ar[r]  & \ker(r) \ar[r] & im(\jmath_+^*) \ar[r]_{r} &
im(\jmath_-^*) \ar[r] & 0 }\end{array}.
\end{equation}
The map $\imath_-^*$ is a surjection, by the inductive hypothesis.
To show that $\imath_+^*$ is a surjection, it suffices to show
that the dotted vertical arrow is a surjection. That $\imath_+^*$
is a surjection then follows by a diagram chase.

We first recall a fact about Euler classes. Suppose $T$ acts on a
complex vector bundle ${\cal E}$ over a manifold $\Sigma$, with fixed point
set precisely $\Sigma$. Decompose ${\cal E}$ into the direct sum of bundles
${\cal E}_\alpha$, where $T$ acts on ${\cal E}_\alpha$ by weight $\alpha\in
\t^*_{\Z}$. Assume that weights $\alpha$ are distinct and pairwise
relatively prime in $H^2_T(pt;\Z)\iso\t^*_{\Z}$. Let $\tau_\alpha$ be the
equivariant Euler class of the subbundle ${\cal E}_\alpha$. Then
if $y \in H^*_T(\Sigma;\Z)$ is a multiple of $\tau_\alpha$ for each
$\alpha$, then $y$ is a multiple of the product of the $\tau_\alpha$.
This follows from the proof of \cite[Lemma~3.2]{TW99}.  Although their
Lemma is stated for $\Q$ coefficients, the argument goes through given
our assumption of relative primality of the weights.

We now characterize $\ker(r)$.  Suppose $\eta$ is a class in
$H_T^*(\overline{N}_{c_i}^+;\Z)$ such that its restriction to
$H_T^*(\Sigma_{c_i}^-;\Z)$ is zero; that is, $r\circ
\jmath^*_+(\eta) = 0$. Let $\nu = \oplus_{\alpha} \nu_{\alpha}$ be
the $T$-invariant decomposition of the negative normal bundle to
$\Sigma_{c_i}$, where the weights $\{\alpha\}$ are distinct. Let
$N_{\alpha}$ be the component of the one-skeleton $N$ such that
the closure contains $\Sigma_{c_i}$ corresponding to the weight
$\alpha$. The closure $\overline{N}_{\alpha}$ is a symplectic
manifold, and the restriction of $\eta$ to
$\overline{N}^+_{\alpha}:= \overline{N}_{\alpha} \cap M_{c_i}^+$
has the property that it vanishes when restricted to
$\overline{N}^-_{\alpha}:=\overline{N}_{\alpha} \cap M_{c_i}^-$.
This is because $\eta$ vanishes when restricted to
$\Sigma_{c_i}^-$ (by injectivity for $\overline{N}_{\alpha}$). By
Proposition~\ref{prop:key} applied to the pair
$(\overline{N}^+_{\alpha}, \overline{N}^-_{\alpha})$, we may
conclude that $\eta$ restricted to $\Sigma_{c_i}$ must be a
multiple of each $\tau_{\alpha}$. By assumption, any two distinct
$T$ weights occurring in the negative normal bundle to
$\Sigma_{c_i}$ are relatively prime in $H^*_T(pt;\Z)$. Hence by
the fact recalled in the previous paragraph, the restriction of
$\eta$ to $H^*_T(\Sigma_{c_i};\Z)$ has to be a multiple of the
product of the $\tau_{\alpha}$, which is the equivariant Euler class
of the negative normal bundle to $\Sigma_{c_i}$.

We now show that the left vertical arrow in the
diagram~\eqref{eq:inductiveStep} is a surjection. We have 
shown that any element in $\ker(r)$ is, when restricted to
$\Sigma_{c_i}$, a multiple of the equivariant Euler class
$\tau_{c_i}$ of the negative normal bundle to $\Sigma_{c_i}$. On the
other hand, any class in $H^*_T(\Sigma_c)$ which is a multiple of
$\tau_c$ is the image of an element in $H^*_T(M^+, M^-)$ by
Proposition~\ref{prop:key}. Hence the left vertical arrow is
surjective, and the surjectivity of $\imath_-^*$ follows by the
five lemma.
\end{proof}

Theorem~\ref{th:oneskeleton} tells us that it suffices to identify
the image of $\jmath^*$ to find a description of $H_T^*(M)$.  We
will now place stronger hypotheses on the fixed point set $\Sigma$
and the one-skeleton $\overline{N}$ so that the image of
$\jmath^*$ has a simple combinatorial description. We make the
following definition.
\begin{definition}\label{def:GKMmanifold}
Let $M$ be a symplectic manifold equipped with a Hamiltonian
$T$ action. We say that the action is {\em GKM} if $M^T$ consists of
finitely many isolated points, and the $T$-isotropy weights
$\alpha_{i,p}$\glossary{alphaip@$\alpha_{i,p}$, $T$-isotropy weights
at a fixed point of $M$} at a given fixed point are pairwise
relatively prime in $H^*_T(M;\Z)$.
\end{definition}
Henceforth, we assume that our action is GKM.  Thus each component of
$\Sigma$ is an isolated point, and all equivariant Euler classes are
elements of $H^*_T(pt;\Z)$, given as products of the relevant isotropy
weights.  Moreover, if the moment map is proper and bounded below in
some direction, the one-skeleton is a union of copies of $\C P^1$ and
$\C$, intersecting in fixed points.  When $M$ is compact, the pairwise
relative primality of the isotropy weights is equivalent to the
one-skeleton being two-dimensional \cite{GZ01}.  The same holds for
GKM actions in the presence of a moment map that is proper and bounded
below in some direction, by a symplectic cutting argument.

We now associate a graph $\Gamma$\glossary{Gamma@$\Gamma$, GKM
graph associated to $M$} to the GKM action on $M$ that encodes the
information necessary to compute the equivariant cohomology of
$M$.  We call this the {\em GKM graph}.  The vertices
$V$\glossary{V@$V$, set of vertices of the GKM graph} of $\Gamma$
are the fixed points $M^T$. The edges $E$\glossary{E@$E$, set of
edges in the GKM graph} of $\Gamma$ correspond to the embedded $\C
P^1$'s.  That is, we include an edge between two fixed points
precisely when they are the two fixed points of a $\C P^1$ in the
one-skeleton.  Each edge $e\in E$ is labeled with the weight
$\alpha_e$\glossary{Alphae@$\alpha_e$, $T$-isotropy weight on $\C
P^1$ corresponding to edge} of the torus action on that copy of
$\C P^1$.  Notice that the $\C$'s in the one-skeleton
equivariantly retract, and therefore do not contribute to the
cohomology of the one-skeleton. Thus, we do not record this
information in the graph $\Gamma$.

The computation of the cohomology of the one-skeleton for a GKM
action now boils down to the computation of the $T$-equivariant
cohomology of $\C P^1$.  For the proof of the following Lemma,
see, for instance, \cite{HHH04}.

\begin{lemma}\label{lem:cp1coh}
Suppose $T$ acts linearly and nontrivially on $\C P^1$ with
weight $\alpha$. Then the inclusion of the fixed points \( (\C P^1)^T=\{N,S\} \into \C
P^1\) induces an injection $ \imath^*: H_T^*(\C P^1;\Z) \to
H_T^*(\{N,S\};\Z), $ with image
\[
\imath^*(H_T^*(\C P^1);\Z) = \left\{ (f, g) \in H_T^*(\{N\};\Z)
\oplus H_T^*(\{S\};\Z)\ \Big| \  (f - g) \cong 0\ (\mathrm{ mod}\
\alpha) \right\}.
\]
\end{lemma}

Motivated by this lemma, we now define the graph cohomology of
$\Gamma$ to be
$$
H^*(\Gamma, \alpha) := \left\{ f:V\to H_T^*(pt;\Z)\ \left| \
\begin{array}{c}f(p)-f(q)\equiv 0\ ({\mathrm{ mod}}\ \alpha_e)\\
\mbox{ for every edge } e = (p,q)
\end{array}\right.
\right\} \subseteq H_T^*( V;\Z).
$$
\glossary{HGamma@$H^*(\Gamma,\alpha)$, the graph cohomology associated
  to $\Gamma$ with edge-weight data $\alpha$}%
Since the one-skeleton consists of $\C P^1$'s (and equivariantly
retractable $\C$'s) intersecting at fixed points, a Mayer-Vietoris
type argument shows that the image of the cohomology of the
one-skeleton under $\jmath^*$ is precisely the graph cohomology.
This, combined with Theorem~\ref{th:oneskeleton}, yields the
following theorem.

\begin{theorem}\label{thm:GKMnoncompact}
Let a torus $T$ act on a symplectic manifold $M$ with moment map
$\mu : M\to \algt^*$ that is proper and bounded below in some
generic direction. Suppose that $\Sigma=M^T$ consists of only
finitely many isolated points, and that the $T$-isotropy weights
at $p\in \Sigma$ are pairwise relatively prime in $H_T^*(pt;\Z)$.
Then under the map $\imath^*$, the equivariant cohomology
$H_T^*(M;\Z)$ maps isomorphically onto $H^*(\Gamma,\alpha)$.
\end{theorem}

\section{Examples: hypertoric varieties}\label{sec:hypertoric}

In this section, we present the examples that motivated the work in
this paper.
These are the {\em hypertoric} varieties studied in
\cite{BieDan, HP02, Kon99, Kon00}. Just as their K\"ahler counterparts,
hypertoric varieties come equipped with natural
$T^d$ actions. However, it is important to note that the GKM
conditions only hold for hypertoric varieties when they are viewed as
$T^d \times S^1$ spaces, where the $S^1$ action is an extra piece of
structure on hypertoric varieties not present in the K\"ahler
versions. This will be explained in detail below. Throughout this
section, we take the coefficient ring \(R = \Z.\)

We first set some notation in order to facilitate discussion of the
examples. For details we refer the reader to \cite{BieDan, HP02}.
Let $T^n$ be the real $n$-dimensional torus acting on $\Cn$, with
induced action on $\Hn\cong\cot$
given by $t(z,w) = (tz, t^{-1} w)$.
Let \(\{a_i\}_{1 \leq i \leq n}\)\glossary{ai@$\{a_i\}$, nonzero
  primitive integer vectors defining the subtorus $T^k \subset T^n$}
be nonzero primitive integer vectors in
\(\liet^d \cong \R^d\) and let $\{\e_i\}$ be the standard basis for
\(\liet^n \cong \R^n\), dual to 
$\{h_i\}$\glossary{hi@$\{h_i\}$, the standard basis for $(\t^n)^*$}
the standard basis for $(\t^n)^*$. Define the map \(\beta: \liet^n \longrightarrow
\liet^d\)\glossary{beta@$\beta$, the linear map \(\t^n \to \t^d\)
  defining the subtorus $T^k$} by setting \(\beta(\e_i)= a_i,\)
This map fits into an exact sequence
\begin{equation}\label{eq:tseq}
0 \longrightarrow \liet^k \stackrel{\iota}{\longrightarrow} \liet^n
\stackrel{\beta}{\longrightarrow} \liet^d\longrightarrow 0,
\end{equation}
where \(\liet^k := {\mathrm ker}(\beta).\)
Exponentiating yields a subtorus $T^k$ of $T^n$.

\begin{figure}[h]
\psfrag{1}{$1$} \psfrag{2}{$2$} \psfrag{3}{$3$} \psfrag{4}{$4$}
\begin{center}
\includegraphics[height=4cm]{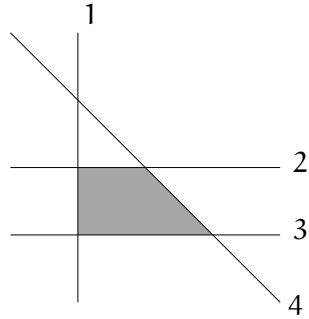}
\end{center}
\begin{center}
\parbox{4.25in}{
\caption{A simple example of a hypertoric variety of real dimension 8
  obtained by reducing $\H^4$ by $T^2$. We label the hyperplane $H_i$
by the index $i$. The region $\Delta$ is
  shaded. The corresponding K\"ahler toric variety is a
  Hirzebruch surface.}\label{fig:ex1}
}
\end{center}
\end{figure}

The of $T^n$ on $\H^n$ is hyperhamiltonian, and so the $T^k$ action is
also hyperhamiltonian. We denote by $M$\glossary{M@$M$, a hypertoric
  variety determined by an affine hyperplane arrangement ${\cal H}$} the hyperk\"ahler reduction
of $\H^n$ by the subtorus $T^k$
at \((\lambda, 0)
\in (\liet^k)^{*}\oplus (\tk_{\C})^*,\) which we assume is a regular value.
This is the hyperk\"ahler analogue of the K\"ahler toric variety $X =
\Cn\mod_{\!\!\lambda}T^k$.\glossary{X@$X$, the K\"ahler toric variety
  associated with the hypertoric variety $M$ and arrangement ${\cal H}$}
The reduction $M$ has a  residual
action of $T^d$ with hyperk\"ahler moment map, denoted \(\mu =
\mu_{\R}\oplus\mu_{\C}.\)
Note that the choice of subtorus $\Tk\subs\Tn$ and a lift
\(\tilde{\lambda}\) of $\lambda$ amounts to choosing an arrangement
${\cal H}$\glossary{H@${\cal H}$, an arrangement of cooriented,
  affine, rational hyperplanes $H_i$} of
cooriented, affine, rational hyperplanes $\{H_i\}_{i=1}^n$, where the $i$th
hyperplane is
\[
H_i = \{x\in\tdd\mid\< x,a_i\>=\< -\tilde{\lambda},\e_i\>\}.
\]
The coorientation comes from knowing for which $x$ we have \(\left<x,
a_i\right> > 0.\) To record the coorientations,
we define the half-spaces
\begin{equation}\label{eq:FGdef}
\begin{array}{c}
F_i\glossary{Fi@$F_i$, half-space defined by $H_i$ containing
    $\Delta$} = \{x\in\tdd\mid\< x,a_i\>\geq\< -\tilde\lambda,\e_i\>\} \\
\text{and}\\
G_i\glossary{Gi@$G_i$, half-space defined by $H_i$ not containing
  $\Delta$}
  = \{x\in\tdd\mid\< x,a_i\>\leq\< -\tilde\lambda,\e_i\>\},
\end{array}
\end{equation}
which intersect in the hyperplane $H_i$.  In our examples, we assume
that the half-spaces $F_i$ intersect in a nonempty bounded polytope
\(\Delta = \cap_{i=1}^n F_i.\) See Figure~\ref{fig:ex1} for an
example. This polytope $\Delta$\glossary{Delta@$\Delta$, polytope
determined by ${\cal H}$} is exactly the image under $\mu_{\R}$ of the
K\"ahler toric variety \(X = \C^n \mod_{\lambda} T^k.\)

  In the case of hypertoric varieties, there is an additional
  residual Hamiltonian $S^1$ action descending from the action of
  $S^1$ on the cotangent bundle $T^{*}\C^n$ that rotates the
  fibers with weight 1. This $S^1$ action is Hamiltonian with respect
  to ${\tilde{\omega}}_{\R}$. Since this restricts to the trivial action on the zero
  section $\C^n$, the $S^1$ action is trivial on the K\"ahler toric
  variety. Hence this action is a new feature of
  hypertoric varieties, and it turns out that this new $S^1$ action is
  essential for the GKM description of the $T^d \times
  S^1$ equivariant cohomology of $M$. We denote the moment map for the extra $S^1$ action by
  $\psi$\glossary{psi@$\psi$, ``extra'' $S^1$ moment map on a
  hypertoric variety}.

  We begin by showing that the hypertoric varieties built above by
  the hyperK\"ahler Delzant construction satisfy the hypotheses of
  Theorem~\ref{thm:GKMnoncompact}. We always assume that the
  affine, rational, cooriented hyperplane arrangement $\arr$ is
  {\em smooth} in the sense of \cite{HP02}, which implies that the
  hypertoric variety $M$ associated to ${\mathcal H}$ is smooth. (In
  particular, this means that the arrangement is {\em simple:} every
  subset of $m$ hyperplanes intersect in codimension $m$.)
  Moreover, we also assume that the polytope $\Delta$ is nonempty and bounded
  in $(\t^d)^*$.  We first show that there is a direction of the
  moment map which is proper, bounded, and Morse.

\begin{lemma}\label{lemma:htmorse}
Let $M$ be the hypertoric variety associated to an affine, cooriented,
rational, smooth hyperplane arrangement ${\cal H}$ such that \(\Delta
= \cap_i F_i\) is nonempty and bounded. Let $(\mu,\psi)$
be the $T^d \times S^1$ moment map on $M$.
Then there is a component of $(\mu,\psi)$ which
is proper, bounded, and Morse.
\end{lemma}

  \begin{proof} Since we assume that $\Delta$ is bounded, by
  \cite[Proposition~1.3]{HP02}, the residual $S^1$ moment map $\psi$ is
  proper. Moreover, since the original $S^1$ moment map $\tilde{\psi}$
  on $T^{*}\C^n$ is given by a norm-square of the cotangent vector, it
  is bounded below by $0$. Hence the moment map $\psi$ on the quotient
  is also bounded below. Now consider the $T^d\times S^1$ moment map
  \((\mu,\psi),\) with values in \((\t^d)^* \times \R \cong
  \R^{d+1}.\) We have just shown that a component $(\mu,\psi)^{\xi} =
  \psi$ of this moment map is proper and bounded below. By taking a
  small enough perturbation of $\xi$, we obtain a generic component of
  the moment map which is proper, bounded, and also Morse.
\end{proof}

We must now show that the \(T^d \times S^1\) fixed points on $M$ is a
finite collection of isolated points, and that the isotropy weights
are relatively prime. We set the
following notation.
The hyperplanes $\{H_i\}$
divide \((\liet^d)^{*} \cong \R^d\) into a finite family of closed,
convex polyhedra
$$\Delta_A = (\cap_{i\notin A}F_i)\cap(\cap_{i\in A}G_i),$$
indexed by subsets $A\subs\{1,\ldots,n\}$\glossary{A@$A$, subset of
  $\{1, 2, \ldots, n\}$ defining a region in $\R^n$}\glossary{DeltaA@$\Delta_A$,
  convex polyhedron in $(\t^d)^*$ associated to \(A \subset \{1,\ldots,n\}.\)}.
For each $A\subs\{1,\ldots,n\}$, let
\[
M_A =
\mr^{-1}(\Delta_A)\cap\mc^{-1}(0).
\]
\glossary{MA@$M_A$, K\"ahler
  subvariety corresponding to the polytope $\Delta_A$}This is a K\"ahler submanifold of
$M$ with respect to $\omega_{\R}$, and is the (possibly noncompact)
K\"ahler toric variety
associated to $\Delta_A$ \cite[6.5]{BieDan}.

\begin{proposition}\label{prop:htGKM}
Let $M$ satisfy the hypotheses of Lemma~\ref{lemma:htmorse}. Then
the action of $T^d\times S^1$ on $M$ is GKM.
\end{proposition}

\begin{proof}
      We will need the following facts, all of which may be found in
      \cite{HP02}.  Since the $\C$ moment map $\mu_{\C}$ is
      $S^1$-equivariant (where $S^1$ acts on $\t_{\C}^*$ by the
      standard rotation action), the $S^1$-fixed points of $M$ must
      lie in \(\mu_{\C}^{-1}(0) = \bigcup_A M_A.\) On each $M_A$, the
      torus $T^d$ acts in a Hamiltonian fashion with respect to $\omega_{\R}$
      with moment map \(\mu_{\R} \vert_{M_A}\) and image
      $\Delta_A$. Moreover, on each $M_A$, the extra $S^1$ action acts
      as a subtorus of $T^d$, determined combinatorially by $A$.

      Since we are looking for $T^d \times S^1$-fixed points, the fact
      that all the $S^1$-fixed points are contained in $\mu_{\C}^{-1}(0)$ allows us
      to restrict our attention to the toric varieties $M_A$. Since
      the $M_A$ are usual toric varieties, we find immediately that
      $M^{T^d \times S^1}$ is a subset of the points in $\mu_{\C}^{-1}(0)$
      corresponding to the vertices \(v \in (\t^d)^{*}\) of the
      polyhedral complex defined by ${\mathcal H}$. On the other hand,
      we know from the description of $M^{S^1}$ in \cite{HP02} that
      each such point in $\mu_{\C}^{-1}(0)$ corresponding to a vertex $v$ is
also fixed by $S^1$. Hence the fixed points $M^{T^d \times S^1}$
      are isolated, with images under $\mu_{\R}$ exactly the vertices
      \(v = \cap_{i \in I} H_i\) in the hyperplane arrangement.  In
      particular, \(|M^{T^d \times S^1}|\) is finite.

\begin{figure}[h]
\begin{center}
\psfrag{1}{$1$}
\psfrag{2}{$2$}
\psfrag{3}{$3$}
\psfrag{4}{$4$}
\includegraphics[height=4cm]{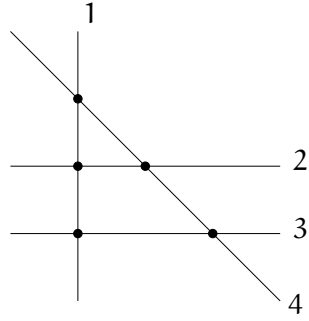}
\end{center}
\begin{center}
\parbox{4.25in}{
\caption{The $T^d \times S^1$-fixed points are mapped to the vertices
of the hyperplane arrangement.}\label{fig:fixedpts}
}
\end{center}
\end{figure}

We must now check that the \(T^d \times S^1\) weights at a given fixed
point $p$ are pairwise relatively prime in $H^*_{T^d\times
S^1}(pt)$. Let \(I \subset \{1, 2, \ldots, n\}\)\glossary{I@$I$, a
subset of $\{1, 2, \ldots n\}$ of size $d$ such that \(\cap_{i \in I}
H_i \neq \emptyset\) } be a subset of size $d$ such that
\(\cap_{i \in I} H_i \neq \emptyset.\) Since ${\cal H}$ is simple, the
intersection is a single vertex $v$. Let $p$ be the fixed point in
$M^{T^d \times S^1}$ corresponding to the vertex \(v = \mu_{\R}(p).\)
We wish to decompose \(T_pM\) under the $T^d \times S^1$-isotropy
action into a sum of 1-dimensional pieces. Since the arrangement
${\mathcal H}$ is simple, there are exactly $2d$ edges coming out of
the vertex $v$, with two edges for each \(i \in I.\) See
Figure~\ref{fig:Tdecomp}. Each edge $e$ defines part of a polytope
$\Delta_A$ corresponding to a subvariety $M_A$ containing $p$. Since
$M_A$ is a standard toric variety, there exists a 1-dimensional weight
space in \(T_pM_A \subseteq T_pM\) with $T^d$ weight $\alpha_e$, where
$\alpha_e$ is the weight corresponding to that edge in
$(\t^d)^{*}$. Since all the weights $\alpha_e$ are distinct in
$(\t^d)^*$, we get a $T^d$ decomposition
\begin{equation}\label{eqn:TpMdecomp}
T_pM \cong \oplus_{i=1}^{2d} \C_{\alpha_{e_i}}.
\end{equation}
This is also a $T^d \times S^1$ decomposition because the $S^1$
commutes with $T^d$.

\begin{figure}[h]
\begin{center}
\psfrag{1}{$1$}
\psfrag{2}{$2$}
\psfrag{3}{$3$}
\psfrag{4}{$4$}
\includegraphics[height=4cm]{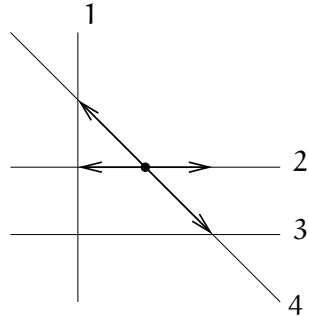}
\end{center}
\begin{center}
\parbox{4.25in}{
\caption{Each edge coming out of a vertex corresponds to a
1-dimensional subspace in $T_pM$. Here, the dimension of the
hypertoric variety $M$ is 4, and there are 4 distinct $T^d$
weights.}\label{fig:Tdecomp}
}
\end{center}
\end{figure}

      We must now show that the $T^d \times S^1$ weights are pairwise
      linearly independent. The hyperplane arrangement $\arr$ is
      simple, so for \(v = \cap_{i \in I} H_i\) as above, the
      collection \(\{a_i\}_{i\in I}\) form a $\Z$ basis of
      \((\t^d)_{\Z}.\) Recall that for each hyperplane $H_i$, we have
      two weights $\alpha_e$ corresponding to $H_i$ in the
      decomposition~\eqref{eqn:TpMdecomp}. These are the two edges
that do not lie in the hyperplane $H_i$. Let \(\{\alpha_{e_i}\}\) be
      a collection of $T^d$ weights in $T_pM$ with \(|\{\alpha_{e_i}\}| =
      d,\) where we have chosen a single weight corresponding to each
      $H_i$. Then the fact that $\arr$ is simple implies that the
      collection of $T^d$ weights \(\{\alpha_{e_i}\}\) is linearly
      independent over $\Z$, so in particular pairwise relatively
      prime over $H^*_{T^d}(pt;\Z)$. We may immediately conclude that
      for $T^d \times S^1$ edge weights \(\alpha_{e_i}, \alpha_{e_j}\)
      (here we abuse notation and use $\alpha_e$ to denote both $T^d$
      and $T^d \times S^1$ weights) are pairwise relatively prime over
      $H^*_{T^d\times S^1}(pt;\Z)$ if $e_i, e_j$ lie on different
      hyperplanes.

      It remains to deal with the case when two weights \(\alpha_{e},
      \alpha_{e\prime}\) correspond to the same hyperplane. In this
      case, as $T^d$ weights, they are negative multiples of one
      another. Hence, to get relative primality, we must compare their
      $S^1$ weights.  In order to compute this $S^1$ weight on a given
      $\C_{\alpha_e}$, we use the fact that the action of $S^1$ on
      each $M_A$ is that of a subtorus (depending combinatorially on
      $A$) of $T^d$. It follows from the computation in \cite{HP02}
      that the $S^1$ weight on $\C_{\alpha_e}$ is given by
      \(\<\alpha_e, - \sum_{i \in A} a_i\> \in \Z\) for $\Delta_A$
      containing both the vertex $v$ and the edge $e$. Although the
      choice of $A$ here is not unique, the weight is
      well-defined. For if $\alpha_e$ is an edge weight for $M_A$ and
      $M_{A'}$, where $A$ and $A'$ differ by a single $i$, then $a_i$
      is necessarily in the annihilator of $\alpha_e$. See
      Figure~\ref{fig:S1weight}. By a simple inductive argument, we
      conclude that the pairing above remains constant for different
      choices of $M_A$.

\begin{figure}[h]
\begin{center}
\psfrag{a}{${\scriptstyle a_j}$}
\psfrag{H}{${\scriptstyle H_j}$}
\psfrag{v}{${\scriptstyle w = - \sum_{i\in A}a_i}$}
\psfrag{w}{${\scriptstyle w' = w - a_j}$}
\psfrag{b}{${\scriptstyle \alpha_e}$}
\includegraphics[height=3cm]{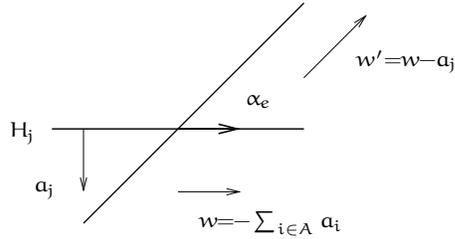}
\end{center}
\begin{center}
\parbox{4.25in}{
\caption{The computation of the $S^1$ weight for the edge $e$. The
      pairing $\<\alpha_e, -\sum_{i\in A} a_i\>$ is well-defined since
      for two adjacent regions, the corresponding vectors $v$ and $v'$
      differ by a vector $a_j$ perpendicular to
      $\alpha_e$.}\label{fig:S1weight}
}
\end{center}
\end{figure}
      To see that $\alpha_e, \alpha_{e'}$ are relatively prime, it
      suffices to check that the $S^1$ weights are not negative
      multiples of each other. Let $A$ be such that $\Delta_A$ contains $v$ and
      $e$. Let $a_j$ define the (unique) hyperplane $H_j, j \in I$,
      for which $\alpha_e, \alpha_{e'}$ do not lie on $H_j$.  Without
      loss of generality, we assume \(\<\alpha_e, a_j\> > 0.\) Then
      \(\<\alpha_{e'}, -\sum_{i \in A} a_i - a_j\> = - \<\alpha_e,
      -\sum_{i \in A} a_i\> + \<\alpha_e, a_j\>.\) See
      Figure~\ref{fig:S1weight2}. Since \(\<\alpha_e, a_j\> \neq 0,\)
      the $S^1$ weights are not negative multiples, and the $T^d\times
      S^1$ weights $\alpha_e, \alpha_{e'}$ are relatively prime.
\begin{figure}[h]
\begin{center}
\psfrag{a}{${\scriptstyle a_j}$}
\psfrag{H}{${\scriptstyle H_j}$}
\psfrag{v}{${\scriptstyle w = - \sum_{i\in A}a_i}$}
\psfrag{w}{${\scriptstyle w' = w - a_j}$}
\psfrag{b}{${\scriptstyle \alpha_e}$}
\psfrag{c}{${\scriptstyle \alpha_{e'}}$}
\includegraphics[height=3cm]{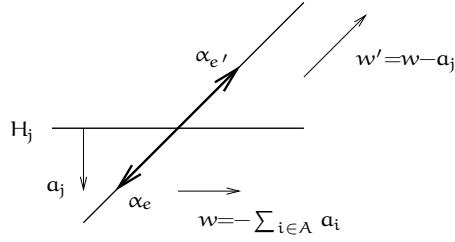}
\end{center}
\begin{center}
\parbox{4.25in}{
\caption{Comparison of the $S^1$ weights for two edges on the same
hyperplane.}\label{fig:S1weight2} }
\end{center}
\end{figure}
\end{proof}

\begin{remark}\label{rmk:notTGKM}
From the proof of Proposition~\ref{prop:htGKM}, it is evident that
the $T^d$ action on $M$ is {\em not} GKM in the sense of
Definition~\ref{def:GKMmanifold}.
\end{remark}

We now give a GKM description of the $T^d \times S^1$-equivariant
cohomology of a hypertoric variety $M$ in the sense of
Section~\ref{sec:gkm}. Let $\Gamma = (V,E)$ denote the GKM graph
of $M$, and let $H^*(\Gamma, \alpha)$ denote its graph cohomology.
By Theorem~\ref{thm:GKMnoncompact}, Lemma~\ref{lemma:htmorse}, and
Proposition~\ref{prop:htGKM}, we may conclude that the image of
the map
\[
\imath^*: H^*_{T^d \times S^1}(M) \longrightarrow H^*_{T^d \times
  S^1}(M^{T^d \times S^1})
\]
induced by inclusion is an injection, with image $H^*(\Gamma,
\alpha)$.

We now have an explicit description of the $T^d\times
S^1$-equivariant cohomology of $M$ as a subring of the sum of
polynomial rings $H^*_{T^d\times S^1}(M^{T^d\times S^1};\Z)$.
Another description, in terms of generators and relations, of the
$T^d\times S^1$-equivariant cohomology of $M$ was given in
\cite{HP02}. We will now give a set of ring generators of
$H^*(\Gamma, \alpha) \cong H^*_{T^d \times S^1}(M;\Z)$ by
constructing an isomorphism between the quotient description of
the \(T^d \times S^1\)-equivariant cohomology given in \cite{HP02}
and the GKM description via $H^*(\Gamma, \alpha)$. We first recall
the following theorem. The $F_i, G_i$ are defined in
equation~\eqref{eq:FGdef}.

\begin{theorem}\label{thm:quotient}{\em\cite[4.4]{HP02}}
Let $M$ be the hypertoric variety satisfying the hypotheses of
Lemma~\ref{lemma:htmorse}.
Given any minimal set $S\subs \{1,\ldots,n\}$\glossary{S@$S$, minimal
  subset of \(\{1,\ldots,n\}\) such that \(\cap_{i \in S} H_i = \emptyset\)} such that
$\cap_{i\in S} H_i = \emptyset$, let $S = S_1\sqcup S_2$
be the unique splitting of $S$ such that
\begin{equation}\label{eq:FGsplit}
\big(\cap_{i\in S_1}G_i\big)\cap
\big(\cap_{j\in S_2}F_j\big) = \emptyset.
\end{equation}
Then the $T^d \times S^1$-equivariant cohomology of $M$ is given by
\begin{equation}\label{eq:ideal}
H^*_{T^d \times S^1}(M;\Z) \cong\Z[u_1,\ldots,u_n,x]\bigg/
\< \prod_{i\in S_1}u_i\times\prod_{j\in S_2}(x-u_j)
\hs\bigg{\vert}\hs\bigcap_{i\in S} H_i = \emptyset\>.
\end{equation}
\end{theorem}

The isomorphism between this quotient description and the GKM
description of $H^*_{T^d\times S^1}(M;\Z)$ which we present below
is similar in spirit to the isomorphism between the corresponding
descriptions for the $T^d$-equivariant cohomology ring of the
K\"ahler toric variety $X$.  The essential geometric insight is to
recognize the generators $u_i$ as the $T^d \times S^1$-equivariant
Chern classes of certain natural line bundles over $M$.

We first set some notation. Let \(v = \cap_{i \in I} H_i\) be a vertex. For each such $v$, we
define the following subsets of $\{1, 2,\ldots, n\}$:
\begin{eqnarray*}
I_v\glossary{Iv@$I_v$, subset of $\{1,\ldots,n\}$} & := & \{i \mid v \in H_i\}, \\
J_v\glossary{Jv@$J_v$, subset of $\{1,\ldots,n\}$} & := & \{i \mid v \in F_i, v \notin H_i\},\\
K_v\glossary{Kv@$K_v$, subset of $\{1,\ldots,n\}$} & := & \{i \mid v \in G_i, v \notin H_i\}.\\
\end{eqnarray*}
Clearly \(I_v = I,\) the three sets $I_v,J_v,K_v$ are
pairwise disjoint, and \(I_v \cup J_v \cup K_v = \{1, 2, \ldots,
n\}.\) For $v$ a vertex and \(i \in I_v,\) we define \(\eta_{v,i} \in
\tdd_{\Z}\)\glossary{etavi@$\eta_{v,i}$, $T^d$-isotropy weight at vertex
  $v$ for $u_i$} to be the element satisfying
\begin{equation}\label{eq:etavi}
\left< \eta_{v,i}, a_j \right> = 0, \quad \forall  j \in I_v\setminus\{i\},
\quad \mathrm{and} \quad \left< \eta_{v,i}, a_i \right> = 1.
\end{equation}
This is well-defined since we assume ${\cal H}$ is simple, so
the vectors $\{a_i\}_{i\in I_v}$
form a $\Z$ basis for $\td_{\Z}$.

We now give a
GKM description of the \(T \times S^1\)-equivariant cohomology of
hypertoric varieties.
We specify a $T^d \times S^1$ weight as a pair \((\alpha, c) \in \tdd_{\Z}
\oplus \Z.\) Let
$x$ denote the equivariantly constant class in $H^*(\Gamma)$
corresponding to the integral basis element for $Lie(S^1)$.

\begin{theorem}\label{thm:quotient_GKM_isom}
Let $M$ be a hypertoric variety satisfying the hypotheses of
Lemma~\ref{lemma:htmorse}, $\mathcal{I}$ the ideal given in \eqref{eq:ideal},
and $H^*(\Gamma)$ denote the graph cohomology associated to $M$.  
Then
the inclusion \(M^{T^d\times S^1} \into M\) induces an isomorphism 
\[
\xymatrix{ H^*_{T^d \times S^1}(M;\Z)
\cong\Z[u_1,\ldots,u_n,x]/\mathcal{I} \ar[r]^(.75){\cong} &
H^*(\Gamma) \\
u_i\glossary{ui@$u_i$, generators of the equivariant cohomology ring
  $H^*_{T^d \times S^1}(M)$} \ar @{|->}[r] &
\rho_i\glossary{rhoi@$\rho_i$, images in $H^*(\Gamma,\alpha)$ of $u_i$ under the
  restriction to the fixed points $M^{T^d \times S^1}$}, \\
x \ar @{|->}[r] & x,
}
\]
where $\rho_i$ is given by
\[
\rho_i(v) =
\begin{cases}
(\eta_{v,i}, \left< \eta_{v,i}, \sum_{j \in K_v} a_j \right> ), &
  \mathrm{if } \hs i \in I_v \\
(0,0), & \mathrm{if } \hs i \in J_v \\
(0, 1), & \mathrm{if } \hs i \in K_v.\\
\end{cases}
\]
\end{theorem}

\begin{proof} Let $\{h_i\}_{i=1}^n$\glossary{hi@$h_i$, standard basis of
  $\tnd_{\Z}$} be the standard basis of $\tnd_{\Z}$, and
  $\tilde{L}_i$\glossary{Litilde@$\tilde{L}_i$, line bundle over
  $T^*\C^n$} be the topologically trivial bundle over $T^*\C^n$ with
  $T^n\times S^1$-equivariant Chern class $h_i$. Let
  $L_i$\glossary{Li@$L_i$, line bundle over $M$} be the quotient
  bundle \(\tilde{L}_i \mid_{\mu_{HK}^{-1}(\alpha,0)} \big/ T^k.\)
  Then the classes $u_i$ are the $T^d\times S^1$-equivariant Chern
  classes of $L_i$ \cite{HP02}. In order to compute the images of
  $u_i$ in \(H^*(\Gamma),\) it suffices to calculate explicitly
  the $T^d\times S^1$ action on each fiber $L_{i,p} := L_i \mid_p$ for
  \(p \in M^{T^d \times S^1}.\)

Let \(v = \cap_{i \in I_v} H_i\) be the vertex corresponding to the
  fixed point $p$. Let \(\pi: Y = \mu_{HK}^{-1}(\alpha,0) \to
  M\)\glossary{Y@$Y$, preimage of $(\alpha,0)$ under
  $\mu_{HK}$}\glossary{pi@$\pi$, projection of
  \(\mu_{HK}^{-1}(\alpha,0)\) to quotient by $T^k$}
  denote the quotient by $T^k$, and let \((z,w) \in Y\) be a preimage
  of the fixed point $p$. By the moment map conditions and by the
  definitions of $I_v, J_v, K_v$, we have
\[
\begin{cases}
z_i = w_i =0 & i \in I_v, \\
z_i \neq 0, w_i = 0 & i \in J_v, \\
z_i = 0, w_i \neq 0 & i \in K_v. \\
\end{cases}
\]
For each \(i \in \{1, 2, \ldots, n\},\) we wish
to compute the restriction of \(u_i = c_1(L_i)\) to the fixed
point $p$ corresponding to the vertex $v$. Let
\(\gamma_{v,i}\)\glossary{gammavj@$\gamma_{v,j}$, $T^d$ weight
component of \(u_i\) restricted to the fixed point $p$
corresponding to $v$} denote the $T^d$ weight component of $u_i
\mid_{p}$. Since the vectors $\{a_j\}_{j \in I_v}$ form a
$\Z$ basis for $\tdd_{\Z}$, in order to completely specify
$\gamma_{v,i}$, it suffices to compute the pairing
$\left<\gamma_{v,i}, a_j\right>$ for all \(j \in I_v.\) Since we
will do our computations on the preimage $\pi^{-1}(p)$, it will be
convenient to do computations with $\beta^*(\gamma_{v,i})$, where
$\beta^*$ is defined by taking the dual of the exact
sequence~\eqref{eq:tseq}. Let \(\{\epsilon_j\}\) denote the
standard basis for $\tn$, and let $t_j$ denote the elements in the
corresponding $S^1$'s in $T^n$.  Let \(((z,w),q)\)\glossary{q@$q$,
coordinate on a line bundle} denote an element in the total space
of the line bundle $\tilde{L}_i$ over the point $(z,w)$.  For
\(j\in I_v, z_j = w_j = 0,\) so the action of $t_j$ on $((z,w),q)$
is given by
\[
t_j \cdot ((z,w),q) =
\begin{cases}
((z,w),q) & j \neq i \\
((z,w),t_i q) & j = i. \\
\end{cases}
\]
Hence the $T^n$ weight $\beta^*(\gamma_{v,i})$ satisfies
\(\left<\beta^*(\gamma_{v,i}),\epsilon_j\right> =
\left<\gamma_{v,i},a_j\right> = 0, \forall j \neq i,
j \in I_v,\) as well as \(\left<\beta^*(\gamma_{v,i}), \epsilon_i\right> =
\left<\gamma_{v,i}, a_i\right> = 1.\) Hence $\gamma_{v,i} =
\eta_{v,i}$, by definition of $\eta_{v,i}$ in~\eqref{eq:etavi}.

We now compute the $S^1$ weight component of $u_i \mid_p$. Recall that
the extra $S^1$ action on $Y \subseteq T^*\C^n$ is given by rotating
the cotangent direction, so for an element $s\in S^1$,
\[
s \cdot (z,w) = (z,sw).
\]
To compute the $S^1$ action on the fiber of $\tilde{L}_i$ over
$(z,w)$, we must find an element in $T^k$ taking $(z,sw)$ back
to $(z,w)$. The subtorus $T^k$ is defined by the exact
sequence~\eqref{eq:tseq}. In particular, an element \(\Lambda = \sum_{j=1}^n c_j
\epsilon_j \in \tk\) if and only if \(\beta(\Lambda) = \sum_{j=1}^n c_j
a_j = 0.\) Observe that \(w_j \neq 0\) exactly when \(j \in K_v\),
and that \(z_j \neq 0\) exactly when \(j \in J_v.\) Hence the
appropriate element in $T^k$ will be an exponential of $\Lambda =
\sum_{j=1}^n c_j \epsilon_j \in \tk$ with the conditions \(c_j =
1\) for \(j \in K_v\) and \(c_j = 0\) for \(j \in J_v.\)
\footnote{We take the convention that the standard action of $S^1$
on $\C$ is given by \(t \cdot z = t^{-1}z.\) See \cite[Appendix
A]{GGK} for an explanation.} Since the $\{a_j\}_{j\in I_v}$ are an
integral basis for $\tdd_{\Z}$, there is a unique integral solution
\(\{m_j\}_{j \in I_v}\) to the equation
\begin{equation}\label{eq:csol}
\sum_{j\in I_v} m_j a_j + \sum_{j\in K_v} a_j = 0.
\end{equation}
The $S^1$ weight on the fiber of $L_i$ is then given by $m_i$ for $i \in I_v$.
Since $\eta_{v,i}$ satisfies the conditions~\eqref{eq:etavi}, the
coefficient $m_i$ can be computed by the pairing \(\left< \eta_{v,i},
-\sum_{j \in K_v} a_j \right>,\) as desired.

Now we take the case \(i \in J_v.\) Observe that $\tilde{L}_i$ has a
$T^n\times S^1$-equivariant section $\tilde{s}_i(z,w) = z_i$, which
descends to a $T^d\times S^1$-equivariant section $s_i$ of $L_i$ with
zero-section $$Z_i := \{[z,w]\in M\mid z_i=0\}.$$ This zero-section has
(real-)moment image $\mr(Z_i) = G_i$. For \(j \in J_v,\) the vertex
$v$ lies in the interior of $F_j$, so the section $s_i$ is nonzero at
$p$. Hence the $T^d\times S^1$ action on the fiber of $L_j$ at $p$ is
trivial, and the $T^d\times S^1$ weight \(u_j \mid_p\) is \((0,0),\)
as desired.

Finally, consider the case \(i \in K_v.\) We first compute the
$T^d$ weight component $\gamma_{v,i}$ of \(u_i \mid_p.\) By the
same argument as for the case \(i \in I_v,\) it suffices to
compute the pairings
\(\left<\beta^*(\gamma_{v,i}),\epsilon_j\right>\) for \(j \in
I_v.\) This time, since \(i \not \in I_v,\) the $T^d$ action on
the fiber is trivial, and
\(\<\beta^*(\gamma_{v,i}),\epsilon_j\right> =
\left<\gamma_{v,i},a_j\right> = 0\) for all \(j \in I_v.\) Hence
the $T^d$ weight component of $u_i \mid_p$ for \(i \in K_v\) is
$0$, as desired. The $S^1$ weight component of $u_i \mid_p$ is
given by $m_i$ for \(i \in K_v\) in the solution~\eqref{eq:csol}.
Hence the $S^1$ weight component is $1$ for \(i\in K_v,\) as
desired.
\end{proof}

\begin{remark}
Interpreted geometrically on the moment map image, Theorem~\ref{thm:quotient_GKM_isom} states
that for a vertex $v$ lying on a hyperplane $H_i$, the $T^d$ weight
component $\eta_{v,i}$ of the restriction \(u_i \mid_{p}\) is
specified by the following conditions:
\begin{enumerate}
\item the $T^d$ weight $\eta_{v,i}$ lies on the edge $\cap_{j \in I_v,
  j\neq i} H_j$; and
\item the $T^d$ weight $\eta_{v,i}$ has positive inner product with the inward-pointing
  normal vector $a_i$, so in particular it always points ``towards''
  $\Delta$.
\end{enumerate}

\end{remark}

\begin{remark}
It is possible to prove that the images of the $u_i$ in $H^*_{T\times
  S^1}(M^{T\times S^1})$
given in Theorem~\ref{thm:quotient_GKM_isom} do indeed satisfy the GKM
conditions. The proof is rather tedious and we do not include it
  here.
\end{remark}

\begin{figure}[h]
\begin{center}
\psfrag{1}{$1$} \psfrag{2}{$2$} \psfrag{3}{$3$} \psfrag{4}{$4$}
\psfrag{( a )}{(a)} \psfrag{( b )}{(b)} \psfrag{( c )}{(c)}
\includegraphics[width=4in]{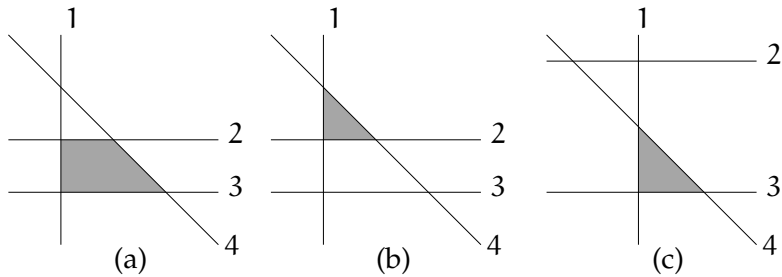}
\end{center}
\begin{center}
\parbox{4.25in}{
\caption{Three different hypertoric varieties.}\label{fig:3exs}
}
\end{center}
\end{figure}

To illustrate Theorem~\ref{thm:quotient_GKM_isom}, we consider the
  hypertoric varieties determined by the hyperplane arrangements in
  Figure~\ref{fig:3exs}.

\begin{example}\label{ex:Ma}
Let $M_a$ denote the hypertoric variety specified by the
hyperplane arrangement in Figure~\ref{fig:3exs}(a). In
\cite{HP02}, the equivariant cohomology $H^*_{T^d\times
S^1}(M_a;\Z)$ is computed to be
\[
H^*_{T^d\times S^1}(M_a;\Z)=\Z[u_1,\ldots,u_4,x]\big/ \< u_2u_3,
u_1(x-u_2)u_4, u_1u_3u_4 \>,
\]
where the $u_i$ are Euler classes of $T^d\times S^1$-equivariant
line bundles $L_i$ over $M_a$.

The images of the $u_i$ are given in Figure~\ref{fig:GKMforMa}. We
choose an integral basis \(\{e_1,e_2\}\) for \(\t^d = \t^2\) as shown in the Figure,
and we denote the integral basis element for $Lie(S^1)$ by $x$. The
equivariantly constant class $x$ maps to the GKM class with weight
$x$ at each point.

\begin{figure}[h]
\begin{center}
\psfrag{w1}{${\scriptstyle e_1}$}
\psfrag{w2}{${\scriptstyle e_2}$}
`\psfrag{w1-w2-w3}{${\scriptstyle e_1-e_2-x}$}
\psfrag{0}{${\scriptstyle 0}$}
\psfrag{w3}{${\scriptstyle x}$}
\psfrag{-w2}{${\scriptstyle -e_2}$}
\psfrag{w1-w2}{${\scriptstyle e_1-e_2}$}
\psfrag{w2}{${\scriptstyle e_2}$}
\psfrag{-w1+w2}{${\scriptstyle -e_1+e_2}$}
\psfrag{-w2-w3}{${\scriptstyle -e_2-x}$}
\psfrag{-w1}{${\scriptstyle -e_1}$}
\psfrag{u1}{$\rho_1$}
\psfrag{u2}{$\rho_2$}
\psfrag{u3}{$\rho_3$}
\psfrag{u4}{$\rho_4$}
\includegraphics[width=5in]{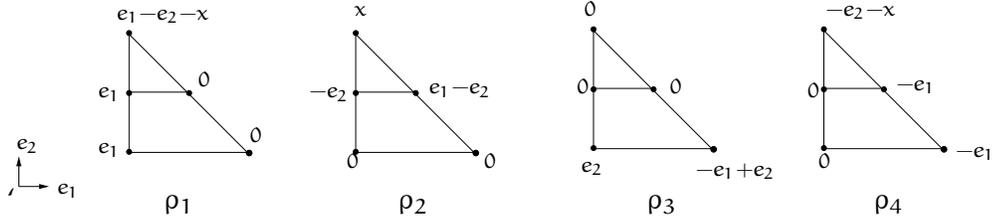}
\end{center}
\begin{center}
\parbox{4.25in}{
\caption{Some ring generators for $H^*_{T^d\times S^1}(M_a;\Z)$.
The plane of the figure is spanned by the two vectors $e_1$ and
$e_2$. The third axis $x$ points out of the
page.}\label{fig:GKMforMa} }
\end{center}
\end{figure}

\end{example}

\begin{example}\label{ex:Mb}

Now let $M_b$ denote the hypertoric variety given by the arrangement
in Figure~\ref{fig:3exs}(b). The equivariant cohomology is computed
\cite{HP02} to be
\[
H^*_{\Td\times S^1}(M_b;\Z)=\Z[u_1,\ldots,u_4,x]\big/ \<
(x-u_2)u_3, u_1u_2u_4, u_1u_3u_4 \>.
\]
We give the GKM descriptions of the ring generators $u_i$ in
Figure~\ref{fig:GKMforMb} below.

\begin{figure}[h]
\begin{center}
\psfrag{w1-w2}{${\scriptstyle e_1-e_2}$}
\psfrag{w1}{${\scriptstyle e_1}$}
\psfrag{0}{${\scriptstyle 0}$}
\psfrag{w2}{${\scriptstyle e_2}$}
\psfrag{-w1+w2}{${\scriptstyle -e_1+e_2}$}
\psfrag{w3}{${\scriptstyle x}$}
\psfrag{w2-w3}{${\scriptstyle e_2-x}$}
\psfrag{-w1+w2-w3}{${\scriptstyle -e_1+e_2-x}$}
\psfrag{-w2}{${\scriptstyle -e_2}$}
\psfrag{-w1}{${\scriptstyle -e_1}$}
\psfrag{u1}{$\rho_1$}
\psfrag{u2}{$\rho_2$}
\psfrag{u3}{$\rho_3$}
\psfrag{u4}{$\rho_4$}
\includegraphics[width=4.5in]{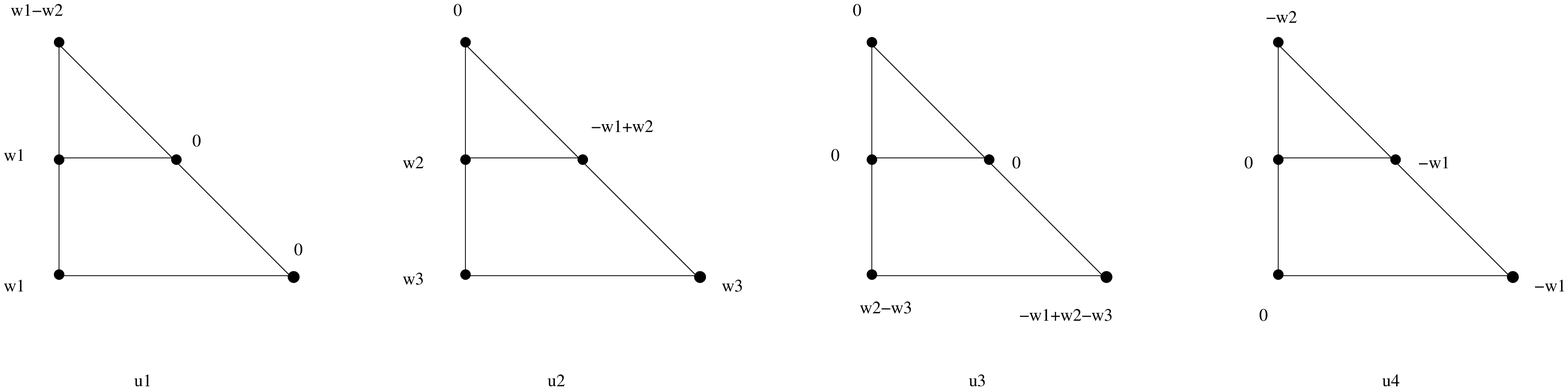}
\end{center}
\begin{center}
\parbox{4.25in}{
\caption{Some ring generators for $H^*_{T^d\times
S^1}(M_b;\Z)$.}\label{fig:GKMforMb} }
\end{center}
\end{figure}

\end{example}

\begin{example}\label{ex:Mc}

Let $M_c$ be the hypertoric variety given by the arrangement
in Figure~\ref{fig:3exs}(c). The equivariant cohomology is computed
\cite{HP02} to be
\[
H^*_{\Td\times S^1}(M_c;\Z)=\Z[u_1,\ldots,u_4,x]\big/ \< u_2u_3,
(x-u_1)u_2(x-u_4), u_1u_3u_4 \>.
\]
We give the GKM image of the $u_i$ in Figure~\ref{fig:GKMforMc} below.

\begin{figure}[h]
\begin{center}
\psfrag{w3}{${\scriptstyle x}$}
\psfrag{w1+w3}{${\scriptstyle e_1+x}$}
\psfrag{w1-w2}{${\scriptstyle e_1-e_2}$}
\psfrag{w1}{${\scriptstyle e_1}$}
\psfrag{0}{${\scriptstyle 0}$}
\psfrag{w1-w2-w3}{${\scriptstyle e_1-e_2-x}$}
\psfrag{-w2-w3}{${\scriptstyle -e_2-x}$}
\psfrag{w2}{${\scriptstyle e_2}$}
\psfrag{-w1+w2}{${\scriptstyle -e_1+e_2}$}
\psfrag{w1-w2+w3}{${\scriptstyle e_1-e_2+x}$}
\psfrag{w3}{${\scriptstyle x}$}
\psfrag{-w2}{${\scriptstyle -e_2}$}
\psfrag{-w1}{${\scriptstyle -e_1}$}
\psfrag{u1}{$\rho_1$}
\psfrag{u2}{$\rho_2$}
\psfrag{u3}{$\rho_3$}
\psfrag{u4}{$\rho_4$}
\includegraphics[width=4.5in]{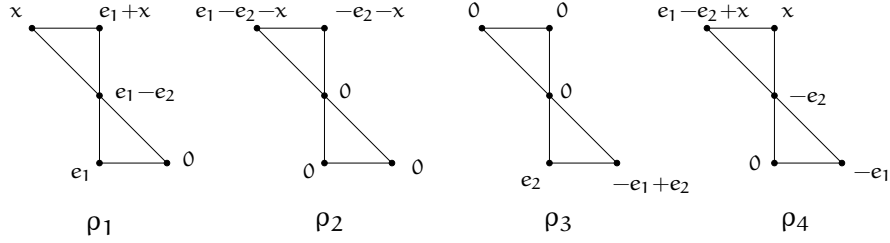}
\end{center}
\begin{center}
\parbox{4.25in}{
\caption{Some ring generators for $H^*_{T^d\times
S^1}(M_c;\Z)$.}\label{fig:GKMforMc} }
\end{center}
\end{figure}

\end{example}

We end this section with a discussion of the $T^d$ action on a
hypertoric variety, considered as a GKM action.  We have already noted
in Remark~\ref{rmk:notTGKM} that the $T^d$ action on $M$ does not
satisfy the GKM hypotheses in the sense of
Definition~\ref{def:GKMmanifold}. It
{\em does} satisfy the more general GKM conditions considered in
\cite{HHH04}, in which the theory is developed in the language of cell
complexes with a compatible $T$ action.  Thus, by \cite[Theorem
3.4]{HHH04}, $H^*_{T^d}(M;\Z)$ does admit a GKM description in
$H^*_{T^d}(M^{T^d};\Z)$. We will now exploit our knowledge of the GKM
description of $H^*_{T^d \times S^1}(M;\Z)$ to give an explicit list
of ring generators for $H^*_{T^d}(M;\Z)$, described as elements of
$H^*_{T^d}(M^{T^d};\Z)$.

\begin{remark} 
Note that the techniques in \cite{HHH04} do not
in general yield ring generators for the $T$-equivariant cohomology, so
this is a new result from our explicit analysis of $M$ as a $T^d\times
S^1$ space.
\end{remark}

We will obtain GKM ring generators for $H^*_{T^d}(M;\Z)$ by ``GKM
in stages.'' First, recall that in
the Borel construction \(M \times_{T^d} ET^d\) for
\(H^*_{T^d}(M;\Z),\) we may use instead of $ET^d$ any contractible
space on which $T^d$ acts freely. In particular, we may use
\(E(T^d \times S^1).\) Hence there is a natural
map
\[
H^*_{T^d\times S^1}(M;\Z) \to H^*_{T^d}(M;\Z)
\]
induced by the inclusion
\[
M\times_{T^d} E(T^d \times S^1) \into M\times_{T^d \times S^1} E(T^d\times S^1).
\]
In our situation, we have in addition that \(M^{T^d} = M^{T^d \times
  S^1},\) so we obtain a commutative diagram
\[
\begin{array}{c}
\xymatrix{ M \times_{T^d \times S^1} E(T^d\times S^1)   &
M^{T^d \times S^1} \times_{T^d \times S^1}
E(T^d \times S^1) \ar@{_{(}->}[l] \\
M \times_{T^d} E(T^d \times S^1) \ar[u] & M^{T^d} \times_{T^d} E(T^d \times S^1)
\ar[u]_{\Pi} \ar@{_{(}->}[l]
 }\end{array},
\]
where $\Pi$ uses the equality \(M^{T^d} = M^{T^d \times S^1}.\) We then have a
diagram on equivariant cohomology
\[
\begin{array}{c}
\xymatrix{
H^*_{T^d\times S^1}(M;\Z) \ar@{^{(}->}[r] \ar[d] & H^*_{T^d\times S^1}(M^{T^d\times S^1};\Z) \ar[d]^{\Pi^*} \\
H^*_{T^d}(M;\Z) \ar@{^{(}->}[r] & H^*_{T^d}(M^{T^d};\Z) \\
}\end{array}.
\]
Since the left vertical arrow is a surjection by formality of
\(H^*_{T^d\times S^1}(M;\Z)\) over $H^*_{S^1}(pt;\Z)$ \cite{HP02},
the right vertical arrow \(\Pi^*\) also gives a surjection on the
images. Moreover, since the images of the $u_i$ generate
\(H^*_{T^d}(M),\) in order to give generators for the GKM
description of $H^*_{T^d}(M;\Z)$ in $H^*_{T^d}(M^{T^d};\Z)$, it
suffices to compute \(\Pi^*(\rho_i),\) where the $\rho_i$ are given
in Theorem~\ref{thm:quotient_GKM_isom}. Note that the map $\Pi^*$
is the map that sends $x$ to $0$. We end the section with an example
of an explicit computation.

\begin{example}\label{ex:Ma_T}

Generators of $H^*_{T^d}(M;\Z)$, considered in
$H^*_{T^d}(M^{T^d};\Z)$, are shown in Figure~\ref{fig:TdGKMforMa}
below. They are the $\pi^*(\rho_i)$ for the $\rho_i$ given in
Example~\ref{ex:Ma}.

\begin{figure}[h]
\begin{center}
\psfrag{w1}{${\scriptstyle e_1}$}
\psfrag{w2}{${\scriptstyle e_2}$}
`\psfrag{w1-w2-w3}{${\scriptstyle e_1-e_2}$}
\psfrag{0}{${\scriptstyle 0}$}
\psfrag{w3}{${\scriptstyle 0}$}
\psfrag{-w2}{${\scriptstyle -e_2}$}
\psfrag{w1-w2}{${\scriptstyle e_1-e_2}$}
\psfrag{w2}{${\scriptstyle e_2}$}
\psfrag{-w1+w2}{${\scriptstyle -e_1+e_2}$}
\psfrag{-w2-w3}{${\scriptstyle -e_2}$}
\psfrag{-w1}{${\scriptstyle -e_1}$}
\psfrag{u1}{$\pi^*(\rho_1)$}
\psfrag{u2}{$\pi^*(\rho_2)$}
\psfrag{u3}{$\pi^*(\rho_3)$}
\psfrag{u4}{$\pi^*(\rho_4)$}
\includegraphics[width=4.5in]{GKMforMa.eps}
\end{center}
\begin{center}
\parbox{4.25in}{
\caption{Ring generators for
$H^*_{T^d}(M_a;\Z)$.}\label{fig:TdGKMforMa} }
\end{center}
\end{figure}

\end{example}

\section{The real locus}\label{sec:mod2gkm}

We now consider the situation in which we have, in addition to a
Hamiltonian $T$-action on $M$, an antisymplectic involution
$\sigma$\glossary{sigmasmall@$\sigma$, an antisymplectic
involution on $M$} on $M$ which anticommutes with the action of
$T$, i.e.
\begin{equation}\label{eq:invol}
\sigma(tx) = t^{-1}\sigma(x), \quad \forall x \in M, \hs \forall t \in T.
\end{equation}
The $\sigma$-fixed points \(Q := M^{\sigma}\)\glossary{Q@$Q$, the real
  locus $M^{\sigma}$ of $M$} in $M$ is a real
$n$-dimensional Lagrangian submanifold of $M$, which
we call the {\em real locus} of $M$. The real locus $Q$ and its
properties (e.g. its image under the $T$ moment map, and its
equivariant cohomology) have been extensively studied; see
\cite{Dui83}, \cite{OS99}, \cite{BGH}, \cite{Sch}. Most of the
known results use the assumption, in addition to certain technical
conditions about the $T$ action, that $M$ (and therefore $Q$) is
compact. We will show in this section that many known results generalize
to the situation in which $M$ is not necessarily compact,
but a component of the moment map is proper and bounded below. Again, our
motivating example is the hypertoric variety with its \(T \times S^1\)
action, which will be discussed in detail in Section \ref{sec:realex}.

We begin our story with an extension of a
theorem which states that the image of
the real locus under the $T$ moment map coincides with that of the
whole manifold $M$, i.e. \(\mu(Q) = \mu(M).\) When the manifold $M$
is compact, this result is due to Duistermaat
\cite{Dui83}.

\begin{proposition}\label{pr:mmap}
Let \((M,\omega)\) be a symplectic manifold with a Hamiltonian
$T$ action, $T$ moment map $\mu$, and $\sigma$ an antisymplectic
involution satisfying~\eqref{eq:invol}. Denote by $Q$ the real locus of
$M$, i.e. \(Q := M^{\sigma}.\) Suppose that there is a component of
the moment map $\mu$ which is proper and bounded below.  Then
\(\mu(Q)=\mu(M).\)
\end{proposition}

\begin{proof}
Let $\mu^{\xi}$ be a component of a moment map for $M$ which is
proper and bounded below. We may assume $\xi$ is rational. Denote by
$S^1_{\xi}$ the subtorus in $T$ generated by $\xi$.
Since $\sigma(tx) =
t^{-1} \sigma(x)$ for all \(x \in M, t \in T ,\) we may assume that
\(\mu^{\xi}(\sigma(x)) = \mu^{\xi}(x),\) for any \(x \in M\)
\cite[2.2]{OS99}.  Without loss of generality we assume $0$ is the
minimum value of $\mu^{\xi}$ on
$M$. Let $\overline{M}_s$ be the symplectic cut space of $M$ at the value
$s>0$ with respect to the action of $S^1_{\xi}$.
$\overline{M}_{c+s}$. Then $\overline{M}_{c+s}$ is equipped
with an antisymplectic involution $\overline{\sigma}_{c+s}$,
descending from the involution
\(\tilde{\sigma}(m,z) = (\sigma(m),\overline{z})\) on \(M \times \C,\) as
well as a Hamiltonian $T$ action, descending from the action of $T$ on
the first factor. These obey the
relation~\eqref{eq:invol}. Denote by $\overline{\mu}$ the $T$ moment map on
$\overline{M}_s$, and let \(\overline{Q}_s :=
(\overline{M}_s)^{\overline{\sigma}_{c+s}}\) be its real locus.
Since $\mu^{\xi}$ is proper, the symplectic
cut space $\overline{M}_s$ is compact. Thus, Duistermaat's theorem
applies, and \(\overline{\mu}(\overline{M}_s) =
\overline{\mu}(\overline{Q}_s).\) On the other hand, the symplectic cut
space $\overline{M}_s$ contains as an open subset the preimage
\((\mu^{\xi})^{-1}((-\infty,c+s)) \subset M\) in the original
manifold, and on this open piece, the involution $\overline{\sigma}_{c+s}$, the
Hamiltonian $T$ action, and the $T$ moment map $\mu$ all agree with those
just defined on $\overline{M}_s$. Since $s$ was arbitrary, we may conclude
that \(\mu(M) = \mu(Q).\)
\end{proof}

\begin{remark}
Note that for the above proposition, we do {\em not} need to assume
that the $T$ action is GKM. We only need that the $T$ moment map is
proper and a component is bounded below.
\end{remark}

We now turn our attention to the mod $2$ GKM theory for the real locus
$Q$. Since the $T$ action on $M$ anticommutes with $\sigma$,
there is a subgroup
\(\{\pm 1\}^n = (\Zt)^n \subseteq T^n,\)
denoted \(T_\R,\)\glossary{TR@$T_{\R}$, ``real'' $d$-dimensional torus
  $(\Zt)^d$} which preserves $Q$.
Thus we can speak of the $T_\R$-equivariant cohomology of $Q$,
and we will show that under certain conditions, we have
an isomorphism of graded rings
\[
\htmzt \cong \htrx
\]
that halves the grading. (For the compact case, see \cite{BGH, Sch}.)

Henceforth we assume that the $T$ action on $M$ is GKM. In order
to get the isomorphism of graded rings described above, 
we will need additional assumptions on the $T$-isotropy weights at
the fixed points. We first set up the notation.
The {\em mod 2 reduction} of a weight
\(\alpha \in \t_{\Z}^{*}\) in the weight lattice of $T$ is defined
to be its image in \(\t_{\Z}^{*}/2 \t_{\Z}^{*}.\) We will denote
by $\overline{\alpha}_{p,i}$ the mod 2 reduction of a $T$ weight
$\alpha_{p,i}$ at a $T$-fixed point $p$.

\begin{definition}
Let $M$ be a manifold equipped with a $T$ action. Then the action
is {\em mod 2 GKM} if it is GKM and, for every \(p \in M^T,\) the
{\em mod 2} reduced weights
\(\{\overline{\alpha}_{p,i}\}_{i=1}^n\) are all distinct and
nonzero.
\end{definition}

\begin{remark}
In \cite{Sch}, the term {\em $\Zt$ pure} is used; this is equivalent
to {\em mod 2 GKM.}
\end{remark}

In Section~\ref{sec:gkm}, we have already shown that the GKM
theorem holds for noncompact GKM actions if certain conditions
hold on the moment map. In order to show the isomorphism of the
two cohomology rings \(\htmzt\) and \(\htrx,\) we need now to show
that the analogous results hold for a noncompact real locus $Q$ in
the case where the $T$ action is also mod 2 GKM. We will use the Morse theory
of the restricted moment map \(\mu|_Q\) on $Q$.  Let $P$ be the
set
$$
P := \{ x\in Q\ | \ \codim_{\Zt}(Stab_{T_\R}(x)) = 1\}.
$$
As before, we define the one-skeleton $\overline{P}$ of the
$T_{\R}$ action on $Q$ to be the closure of $P$
\glossary{P@$\overline{P}$, the one-skeleton of the real locus
$Q$}.

We first claim that when the $T$ action is mod $2$ GKM, then 
\begin{equation}\label{eq:fixedpts}
M^T = Q^{T_\R}
\end{equation}
and
\begin{equation}\label{eq:one-skel}
\overline{P} = \overline{N}\cap Q.
\end{equation}
Thus the combinatorics of the mod 2 one-skeleton $\overline{P}$ for the $T_{\R}$ action
on $Q$ is the same as that of the one-skeleton $\overline{N}$. We give below a
sketch of a proof of the equality~\eqref{eq:fixedpts} because there is a gap in
its proof in~\cite[Theorem~5.2]{BGH}.  We outline the argument in
\cite[Proposition~5.1.6]{Sch}, and we include this here  because
Schmid's thesis is not available in print.

\begin{proposition}
Let a torus $T$ act on a symplectic manifold $M$ with moment map $\mu
: M\to \algt^*$ that is proper and bounded below in some generic
direction. Suppose further that $M$ is equipped with an antisymplectic
involution $\sigma$ that anticommutes with the $T$ action, and that
the $T$ action is mod $2$ GKM.  Then $M^T = Q^{T_\R}$.
\end{proposition}

\begin{proof}
We must first show that $Q^{T_\R}\subseteq M^T$.  This is what is
shown in the proof of \cite[Theorem~5.2]{BGH}.  The proof of this uses
the fact that the isotropy weights have nonzero mod $2$
reductions. The gap in the proof is that the reverse inclusion $M^T
\subseteq Q^{T_\R}$ is not addressed.  We complete that now.

It suffices to show that $M^T \subseteq Q$.  We first show that there
is at least one $T$-fixed point in $Q$.  Let $\mu^\xi$ be a generic
direction of $\mu$ that is proper and bounded below. Since the fixed
points are isolated, there is exactly one fixed point $p\in M^T$
mapping to the minimum value of $\mu^\xi$.  By
Proposition~\ref{pr:mmap}, we have $\mu^\xi(Q)=\mu^\xi(M)$, and so we
must have $p\in Q$.

Now we show that every $T$-fixed point is in $Q$.  Let $(p,q)$ be an
edge in $\Gamma$, corresponding to an embedded $\C P^1$, where $p$ is
a vertex known to be in $Q$.  This $\C
P^1$ is fixed by a codimension $1$ subtorus $T'\leq T$.  It is a
connected component of $M^{T'}$, and $M^{T'}$ is preserved by
$\sigma$.  As $p\in Q$, this copy of $\C P^1$ must itself be preserved
by $\sigma$.  Applying Proposition~\ref{pr:mmap} to this $\C P^1$
allows us to conclude that $q$ is also an element of $Q$.  Finally,
because the one-skeleton is connected, it follows that every $T$-fixed
point is in $Q$.  This completes the proof.
\end{proof}

\begin{remark}
In the compact setting, \eqref{eq:one-skel}
is proved in \cite[Theorem~5.2]{BGH} and in
\cite[Proposition~5.1.5]{Sch}.  In the noncompact setting, 
this follows by a cutting argument similar to that given in the
proof of Proposition~\ref{pr:mmap}.
\end{remark}

We assume from now on that the $T$ action on our manifold $M$ is 
mod $2$ GKM. In this situation, the moment map $\mu$,
restricted to $Q$, also behaves quite nicely. By our assumptions on
$M$, there is a component of the moment map $\mu^{\xi}$ which is
proper and bounded below, and (since the action is GKM) is a Morse function
on $M$.  In \cite{Dui83}, Duistermaat showed that the restriction of
$\mu^{\xi}$ to $Q$ is also a Morse function on $Q$, and has critical
points exactly \(M^T\cap Q\) when $M$ is compact.  In our situation,
his argument goes through to show that components of $\mu$ are again
Morse functions on $Q$.  Thus, this allows us to compute the equivariant
cohomology $H^*_{T_\R}(Q;\Zt)$ using an argument very similar to that
given in Section~\ref{sec:gkm}.

The statement of the analogous mod $2$ GKM theorem for real loci will
require the definition of a mod $2$ GKM graph \(\Gamma_{\R}\).
The vertices $V_{\R}$ of $\Gamma_{\R}$ are the fixed points
$Q^{T_{\R}}$. The edges are given by the components of the
one-skeleton
$P$ the closure of which is an $S^1$, and they connect the two fixed
points in $Q^{T_{\R}}$ contained in this $S^1$. To each edge $e$ we
associate a weight $\overline{\alpha}_e$ of the $T^{\R}$ action on
this $S^1$. We then define the graph
cohomology $H^*(\Gamma_{\R},\alpha_{\R})$ to be
\begin{eqnarray*}
H^*(\Gamma_{\R}, \alpha_{\R})  & = & \left\{ f:V\to H_{T_{\R}}^*(pt;\Zt)\ \left| \
\begin{array}{c}
f(p)-f(q)\equiv 0\ ({\mathrm{ mod}}\ \overline{\alpha}_e)\\
\mbox{ for every edge } e = (p,q)
\end{array}
\right.\right\}\\
& \subseteq &  H_T^*( V;\Zt).
\end{eqnarray*}
\glossary{HGammaR@$H^*(\Gamma_{\R},\alpha_{\R})$, the mod $2$ graph cohomology associated
  to $\Gamma_{\R}$ with edge-weight data $\overline{\alpha}$}
Note that for \(T_{\R} = (\Zt)^d,\) the $T_{\R}$-equivariant cohomology of a point
$H^*_{T_{\R}}(pt;\Zt)$ is a polynomial ring over $\Zt$ with $d$
generators, where the generators are degree $1$ instead of degree $2$.

We will need the following mod $2$ version of the Atiyah-Bott lemma.
For a proof, see \cite[Proposition~5.3.7]{AllPup} or
\cite[Lemma~2.3]{GolHol04}. 

\begin{lemma}[Atiyah-Bott]\label{lem:ABmod2}
Let ${\cal E}\to B$ be a real rank $\ell$ vector bundle over a compact manifold
$B$. Let $T_{\R}$
be the group $T_{\R} = (\Zt)^d$.
Suppose that $T_{\R}$ acts on ${\cal E}$ with fixed point set
precisely $B$.  Choose a $T_{\R}$-invariant Riemannian metric on
${\cal E}$, and let $D$ and $S$ be the corresponding disk and sphere bundles,
respectively, of ${\cal E}$.  Then the long exact sequence of the pair
$(D,S)$ splits into short exact sequences
$$
\xymatrix{ 0 \ar[r] &  H_{T_{\R}}^*(D,S; \Zt) \ar[r] &  H_{T_{\R}}^*(D;\Zt) \ar[r] &
H_{T_{\R}}^*(S;\Zt) \ar[r] & 0 }.
$$
\end{lemma}

We now prove the real locus version of Theorem~\ref{thm:GKMnoncompact}.

\begin{theorem}\label{thm:mod2gkm}
Let $M$ be a Hamiltonian $T$ space with moment map $\mu$. Assume
that there is a generic direction $\xi$ such that $\mu^\xi$ is
proper and bounded below.  Suppose further that $M$ is equipped
with an antisymplectic involution $\sigma$ that anticommutes with the
$T$ action, and that the $T$ action is mod 2 GKM. Let $Q$ denote the real
locus.  Then the inclusion $Q^{T_{\R}}\into Q$ induces an injection
$$
H_{T_{\R}}^*(Q;\Zt)\into H_{T_{\R}}^*(Q^{T_{\R}};\Zt),
$$
in equivariant cohomology, and the image is precisely
$H^*(\Gamma_\R,\overline{\alpha})$.
\end{theorem}

\begin{proof}
The outline of the proof is the same as that given in Section~\ref{sec:gkm}. We will
only mention the relevant steps where some additional argument is necessary.

We begin with the mod 2 version of Proposition~\ref{prop:key}. To get
the statement, we replace $M$ with the real locus $Q$, $f$ with \(g := f
|_Q,\) $T$ with $T_{\R}$, and Euler classes with Stiefel-Whitney
classes. This follows from Lemma~\ref{lem:ABmod2} by the same argument
as in Proposition~\ref{prop:key}.  Thus the restriction from
\(H_{T_{\R}}^*(Q_c^+;\Zt)\) to \(H_{T_{\R}}^*(Q_c^-;\Z)\) induces an
isomorphism from the kernel of $k^*$ to those classes in
$H_{T_{\R}}^*(p;\Zt)$ which are multiples of $\kappa_p$, the
equivariant Stiefel-Whitney class of the negative normal bundle to $p$ with
respect to $g$.  We are assuming the action is mod $2$ GKM, so the fixed
point components are isolated points.  The mod $2$ version of the
injectivity theorem, Theorem~\ref{thm:injectivity}, now follows by the
same argument, using the $T_{\R}$-equivariant Thom isomorphism theorem
with $\Zt$ coefficients.

Before proceeding to the mod 2 analogue of
Theorem~\ref{th:oneskeleton}, we first take a moment to analyze the
$T_{\R}$-isotropy weights at the fixed points $Q^{T_{\R}}$, which are
the critical points of $g$.
Let $p$ be a fixed point in $M^T = Q^{T_{\R}}$.
There exists a neighborhood of $p$ equivariantly symplectomorphic
to a neighborhood of $0$ in $T_{p}M$ with the symplectic form
$\omega_p$.  Moreover, since \(p \in Q,\) the involution $\sigma$
acts on $T_{p}M$, anticommuting with the action of $T$. The local
normal form theorem in \cite[Theorem~7.1]{OS99} implies that
there exists a $T$-invariant, $\sigma$-antiinvariant compatible
complex structure on $T_{p}M$ making it a complex vector space,
and as a $T_\R$ module, $T_{p}M$ is canonically isomorphic to
the complexification of $(T_{p}M)^{\sigma}$.
More specifically, we have local coordinates such that a neighborhood
of $p$ in $M$ is of the form
\(\oplus_{i=1}^n
\C_{{\alpha}_{i,p}},\) where the $\alpha_{i,p}$ are the
$T$-isotropy weights at $p$, and
$\sigma$ is given by complex
conjugation on each factor. Then a neighborhood of $p$ in $Q$ in these
coordinates is of the form \(\oplus_{i=1}^n
\R_{\overline{\alpha}_{i,p}}.\) Thus the $T_{\R}$-isotropy weights at
$p$ of $T_pQ$ are exactly the mod 2 reductions of the $T$-isotropy weights
$\{\alpha_{i,p}\}$.
In particular, the
$T_{\R}$-equivariant Stiefel-Whitney class of the negative normal bundle in
$Q$ with respect to $g$ is given by the product of the $T_{\R}$ weights
$\overline{\alpha}_{i,p}$ for which \(\<{\alpha}_{i,p}, \xi\> < 0.\)
Since by assumption, all $T_{\R}$ weights are nonzero, the
product is also nonzero in $H^{*}_{T_{\R}}(p; \Zt)$, and
therefore not a zero divisor. Finally, we note that two elements
\(\overline{\alpha}, \overline{\alpha'}\) in
$H^*_{T_{\R}}(p;\Zt)$ are relatively prime if they are nonzero and
distinct.

Using the above observations, the mod $2$ version of
Theorem~\ref{th:oneskeleton} follows from the same argument as in
Section~\ref{sec:gkm}.
\end{proof}

In order to compute $H^*_{T_{\R}}(Q;\Zt)$, it now suffices to compute
the image of $H^*_{T_{\R}}(\overline{P};\Zt)$ in the ring
$H^*_{T_{\R}}(Q^{T_{\R}};\Zt)$. Finally, to observe the isomorphism
between the graded rings \(H^{2*}_T(M;\Zt)\) and
\(H^*_{T_{\R}}(Q;\Zt),\) it suffices to compare the relevant
graphs. Note that since \(M^T = Q^{T_{\R}},\) the vertices of the
graphs $\Gamma$ and $\Gamma_{\R}$ are the same, and because
$\overline{P} = \overline{N}\cap Q$, the graphs are the same.
Moreover, by the argument given in the proof of
Theorem~\ref{thm:mod2gkm}, the images mod 2 of the isotropy
weights $\alpha_e$ on the edges of $\Gamma$ are exactly the
$T_{\R}$-isotropy weights on the edges of $\Gamma_{\R}$. Thus the
combinatorial data specified by the GKM and mod 2 GKM graphs are identical,
and the following corollary is immediate.

\begin{corollary}\label{cor:GKMisom}
Let $M$ be a Hamiltonian $T$ space with moment map $\mu$. Assume
that there is a generic direction $\xi$ such that $\mu^\xi$ is
proper and bounded.  Suppose further that $M$ is equipped with an
antisymplectic involution $\sigma$ that anticommutes with the $T$
action, and that the $T$ action is mod $2$ GKM.
Let $Q$ denote the real locus. Then there is an isomorphism
$$
H_T^{2*}(M;\Zt)\iso H_{T_{\R}}^*(Q;\Zt)
$$
that halves degrees.
\end{corollary}

\section{Examples: real loci of hypertoric varieties}\label{sec:realex}

The hypertoric varieties in Section~\ref{sec:hypertoric} have
a natural antisymplectic involution $\sigma$, induced from the
antisymplectic involution on $T^*\C^n$ given by \((z,w) \mapsto
(\overline{z}, \overline{w}).\) We now analyze the topology of the real
locus of $M$ using techniques of the previous section.

Let $M$ be a hypertoric variety specified by a hyperplane
arrangement ${\cal H}$ and parameter $\alpha$. Then the real locus
$Q$ of $M$ is the set
\[
Q := M^{\sigma} = \{[z,w] \in M \mid z,w \in \R\}.
\]
Since the hypertoric variety has an action of $T^d \times S^1$, the
group acting on $Q$ is now $T_\R = T^d_{\R} \times
\Zt$.

In this situation, we claim that the isomorphism between the
$T_\R$-equivariant cohomology of $Q$ and the $T \times
S^1$-equivariant cohomology of $M$ (both with $\Zt$ coefficients)
can be explicitly described in terms of the line bundles $L_i$
over $M$.

\begin{proposition}\label{prop:htRealLocusIsom}
Let $M$ be a hypertoric variety specified by ${\cal H}$ and
$\alpha$. There are antisymplectic involutions $\sigma_i$
on the total
spaces of the $L_i$, extending the natural involution $\sigma$ on $Q$,
so that the fixed point sets $L_i^\sigma$ are
real vector bundles over $Q$. Moreover, under the isomorphism in
Corollary~\ref{cor:GKMisom}, the Chern class $u_i$ of $L_i$ is
mapped to the Stiefel-Whitney class $\kappa_i$ of the real bundle
$L_i^\sigma$.
\end{proposition}

\begin{proof}
For each line bundle $\tilde{L}_i$ over $T^*\C^n$, define an
involution $\tilde{\sigma}_i$ in coordinates by
$$
\tilde{\sigma}_i\cdot (z,w,q) := (\overline{z} , \overline{w},
\overline{q} ).
$$
This is a lift of the standard antisymplectic involution on
$T^*\C^n$. The $\tilde{\sigma}_i$-fixed point set in $T^*\C^n$ is a
$T_{\R}$-equivariant real line bundle over $T^*\R^n$, and its
complexification is the restriction of $\tilde{L}_i$ to $T^*\R^n$.
Since the $\tilde{\sigma}_i$ are anti-$T$-equivariant, they descend to
antisymplectic involutions $\sigma_i$ on the $L_i$ on
$M$. The $\sigma_i$-fixed point sets are
$T_{\R}$-equivariant real line bundles over the real locus $Q$, and
their complexifications are $L_i |_Q$.

Since the complexification of $(L_i)^{\sigma_i}$ is isomorphic as a
real bundle to $L_i^{\sigma_i} \oplus L_i^{\sigma_i}$, we have
$
\kappa_2(L_i^{{\sigma}_i}) = \kappa_1(L_i^{{\sigma}_i})^2.
$
Under the
natural homomorphism \(H^2(Q;\Z) \to H^2(Q;\Zt)\), the image of the
Chern class of a complex line bundle is the second 
Stiefel-Whitney class, so we may conclude that the isomorphism
between \(H^*_T(M;\Z)\) and \(H^*_{T_\R}(Q;\Zt)\) takes the mod
$2$ Chern class of \(L_i |_Q\) to the Stiefel-Whitney class
\(\kappa_1\) of \(L_i^{\sigma_i}.\)
\end{proof}

\begin{remark}
The presentation of the $T^d \times S^1$-equivariant cohomologies
given in Examples~\ref{ex:Ma},~\ref{ex:Mb},
and~\ref{ex:Mc} is therefore identical to that of the $T_{\R}^d \times
\Zt$-equivariant cohomologies of their real loci, where we use $\Zt$
coefficients and divide all degrees of the classes in half.
\end{remark}

\begin{remark}
The techniques developed in Section~\ref{sec:mod2gkm} and
this description of $H^*_{T_{\R}^d \times \Zt}(Q;\Zt)$ is used in
\cite{HP02} to compute a deformation of the Orlik-Solomon algebra of a
smooth real hyperplane arrangement, depending nontrivially on the
affine structure of the arrangement.
\end{remark}


%
%


\end{spacing}

\end{document}